\documentclass[12pt,a4paper]{article}
\usepackage{amsmath,amsfonts,amssymb,amsthm}
\usepackage{fullpage}

\title{The Dirac operator on $SU_q(2)$}

\author{Ludwik D\c{a}browski,$^1$ \
Giovanni Landi,$^2$ \ 
Andrzej Sitarz,$^3$\footnotemark \\
Walter van Suijlekom,$^1$ \ 
Joseph C. V\'arilly$^4$\footnotemark \\[12pt]
$^1\,$Scuola Internazionale Superiore di Studi Avanzati,\\
      Via Beirut 2-4, 34014 Trieste, Italy\\[3pt]
$^2\,$Dipartimento di Matematica e Informatica,
      Universit\`a di Trieste,\\
      Via Valerio 12/b, 34127 Trieste\\
and INFN, Sezione di Napoli, Napoli, Italy\\[3pt]
$^3\,$Institute of Physics, Jagiellonian University,\\
      Reymonta 4, 30-059 Krak\'ow, Poland\\[3pt]
$^4\,$Departamento de Matem\'atica,
      Universidad de Costa Rica,\\
      2060 San Jos\'e, Costa Rica}

\date{23 January 2005}

\makeatletter
\def\section{\@startsection{section}{1}{\z@}{-3.5ex plus -1ex minus
  -.2ex}{2.3ex plus .2ex}{\large\bf}}
\def\subsection{\@startsection{subsection}{2}{\z@}{-3.25ex plus -1ex
  minus -.2ex}{1.5ex plus .2ex}{\normalsize\bf}}
\makeatother

\newtheorem{thm}{Theorem}[section]
\newtheorem{prop}[thm]{Proposition}
\newtheorem{lem}[thm]{Lemma}

\theoremstyle{definition}
\newtheorem{defn}[thm]{Definition}
\newtheorem{rem}[thm]{Remark}

\numberwithin{equation}{section}

\newcommand{\A}{\mathcal{A}}        
\renewcommand{\a}{\alpha}           
\newcommand{\aappr}{\underline{\alpha}} 
\newcommand{\B}{\mathcal{B}}        
\renewcommand{\b}{\beta}            
\newcommand{\bappr}{\underline{\beta}} 
\newcommand{\braCket}[3]{\langle#1\mathbin|#2\mathbin|#3\rangle}
\newcommand{\C}{\mathbb{C}}         
\newcommand{\CGq}[6]{C_q\!\begin{pmatrix}#1&#2&#3\\#4&#5&#6\end{pmatrix}}
\newcommand{\co}[2]{#1_{(#2)}}      
\newcommand{\cop}{\Delta}           
\newcommand{\dn}{{\mathord{\downarrow}}} 
\DeclareMathOperator{\Dom}{Dom}     
\newcommand{\Dslash}{{D\mkern-11.5mu/\,}} 
\newcommand{\eps}{\varepsilon}      
\renewcommand{\H}{\mathcal{H}}      
\newcommand{\half}{{\mathchoice{\oh}{\oh}{\shalf}{\shalf}}} 
\DeclareMathOperator{\id}{id}       
\newcommand{\K}{\mathcal{K}}        
\newcommand{\ket}[1]{|#1\rangle}    
\newcommand{\kett}[1]{|#1\rangle\!\rangle} 
\newcommand{\la}{\lambda}           
\newcommand{\lt}{\triangleright}    
\newcommand{\N}{\mathbb{N}}         
\newcommand{\nn}{\nonumber}         
\newcommand{\oh}{{\tfrac{1}{2}}}    
\newcommand{\Onda}{\widetilde}      
\newcommand{\ooh}{{\tfrac{3}{2}}}   
\newcommand{\otto}{\leftrightarrow} 
\newcommand{\ox}{\otimes}           
\newcommand{\piappr}{\underline{\pi}^\prime} 
\newcommand{\piop}{{\pi^\circ}}     
\newcommand{\piiapprop}{\underline{\pi}^{\prime\circ}} 
\newcommand{\piiop}{{\pi^{\prime\circ}}} 
\newcommand{\rt}{\triangleleft}     
\newcommand{\sepword}[1]{\quad\mbox{#1}\quad} 
\newcommand{\set}[1]{\{\,#1\,\}}     
\newcommand{\Sf}{\mathbb{S}}        
\newcommand{\sg}{\sigma}            
\newcommand{\shalf}{{\scriptstyle\frac{1}{2}}} 
\DeclareMathOperator{\Spin}{Spin}   
\newcommand{\sq}{\unskip\nobreak\kern5pt\nobreak\vrule
     height4pt width4pt depth0pt}   
\newcommand{\ssesq}{{\scriptstyle\frac{3}{2}}} 
\newcommand{\sesq}{{\mathchoice{\ooh}{\ooh}{\ssesq}{\ssesq}}} 
\DeclareMathOperator{\tsum}{{\textstyle\sum}} 
\newcommand{\U}{\mathcal{U}}        
\newcommand{\up}{{\mathord{\uparrow}}} 
\newcommand{\updn}{{\mathord{\updownarrow}}} 
\newcommand{\vth}{\vartheta}        
\newcommand{\x}{\times}             
\newcommand{\Z}{\mathbb{Z}}         
\newcommand{\7}{\dagger}            
\renewcommand{\.}{\cdot}            
\renewcommand{\:}{\colon}           
\def\<#1,#2>{\langle#1,#2\rangle}   

\begin{document}

\maketitle

\thispagestyle{empty}

\begin{abstract}
We construct a $3^+$-summable spectral triple $(\A(SU_q(2)),\H,D)$
over the quantum group $SU_q(2)$ which is equivariant with respect to
a left and a right action of $\U_q(su(2))$. The geometry is
isospectral to the classical case since the spectrum of the operator
$D$ is the same as that of the usual Dirac operator on the
$3$-dimensional round sphere. The presence of an equivariant real
structure $J$ demands a modification in the axiomatic framework of
spectral geometry, whereby the commutant and first-order properties
need be satisfied only modulo infinitesimals of arbitrary high order.
\end{abstract}

\vspace{2pc}



\renewcommand{\thefootnote}{\fnsymbol{footnote}}
\addtocounter{footnote}{1}
\footnotetext{Partially supported by Polish State
Committee for Scientific Research (KBN) under grant 2\,P03B\,022\,25.}
\addtocounter{footnote}{1}
\footnotetext{Regular Associate of the Abdus Salam ICTP, Trieste.}

\newpage
\section{Introduction}
\label{sec:intro}

In this paper, we show how to successfully construct a
(noncommutative) $3$-dimensional spectral geometry on the manifold of
the quantum group $SU_q(2)$. This is done by building a $3^+$-summable
spectral triple $(\A(SU_q(2)),\H,D)$ which is equivariant with respect
to a left and a right action of $\U_q(su(2))$. The geometry is
isospectral to the classical case in the sense that the spectrum of
the operator $D$ is the same as that of the usual Dirac operator on
the $3$-sphere $\Sf^3 \simeq SU(2)$, with the ``round'' metric.

The possibility of such an isospectral deformation was suggested in
\cite{ConnesL} were the operator $D$ was named the ``true Dirac'' operator. Subsequent  investigations \cite{Goswami} seemed to
rule out this deformation because some of the commutators $[D,x]$,
with $x \in \A(SU_q(2))$, failed to extend to bounded operators, a
property which is essential to the definition of a spectral
triple~\cite{Book}.

These difficulties are overcome here by constructing on a Hilbert
space of spinors $\H$ a spin representation of the algebra
$\A(SU_q(2))$ which differs slightly from the one used in
\cite{Goswami}. Our spin representation is determined by requiring
that it be equivariant with respect to a left and a right action of
$\U_q(su(2))$, a condition which is not present in the previous
approach. The role of Hopf-algebraic equivariance in producing
interesting spectral triples has already met with some success
\cite{ChakrabortyPEqvt,DabrowskiSStwoq}; for a programmatic
viewpoint, see~\cite{SitarzEqvt}.

Our construction of an isospectral noncommutative geometry on the
manifold of $SU_q(2)$, which deforms the usual geometry on the
$3$-dimensional sphere, belongs to an interesting terrain where
noncommutative geometry meets the underlying ``spaces'' of quantum
groups. Recent examples
\cite{DabrowskiLPS,DabrowskiSStwoq,NeshveyevT,SchmuedgenW} are concerned with the
``two-dimensional'' spheres of Podle\'s~\cite{Podles} and more general
flag manifolds \cite{Kraehmer}. The left-equivariant
spectral triple on $SU_q(2)$ constructed in \cite{ChakrabortyPEqvt}
and fully analyzed in \cite{ConnesSUq} is not isospectral and does not
have a good limit at the classical value of the deformation parameter.

\vspace{6pt}

After a brief review in Section~\ref{sec:alg-defns} of $SU_q(2)$ and
its symmetries, mainly to fix notation, we construct its left regular
representation in Section~\ref{sec:eqvt-repns} via equivariance, and
transfer that construction to spinors in
Section~\ref{sec:spinor-repn}. On the Hilbert space of spinors, we
consider in Section~\ref{sec:Dirac-oper} a class of equivariant
``Dirac'' operators $D$. For such an operator $D$ having a classical
spectrum, that is, with eigenvalues depending linearly on ``total
angular momentum'', we prove boundedness of the commutators $[D,x]$,
for all $x \in \A(SU_q(2))$. In fact, this equivariant Dirac operator
is essentially determined by a modified first-order condition, as is
shown later on.

Since the spectrum is classical, the deformation --from $SU(2)$ to
$SU_q(2)$-- is isospectral, and in particular the metric dimension of
the spectral geometry is~$3$.

The new feature of the spin geometry of $SU_q(2)$ is the nature of the
real structure $J$, whose existence is addressed in
Section~\ref{sec:J-oper}. An equivariant $J$ is constructed by
suitably lifting to the Hilbert space of spinors $\H$ the antiunitary
Tomita conjugation operator for the left regular representation of
$\A(SU_q(2))$. However, this $J$ is not the Tomita operator for the
spin representation; for if it were, the spectral triple would inherit
equivariance under the co-opposite symmetry algebra $\U_{1/q}(su(2))$,
forcing it to be trivial. Therefore, the equivariant $J$ we shall use
does not intertwine the spin representation of $\A(SU_q(2))$ with its
commutant, and it is not possible to satisfy all the desirable
properties of a real spectral triple as set forth in
\cite{ConnesReal,Polaris}. This rupture was already observed in
\cite{DabrowskiLPS}; just as in that paper, we must also weaken the
first-order requirement on~$D$.

In Section~\ref{sec:spec-tri}, we rescue the formalism by showing that
the commutant and first-order properties nevertheless do hold, 
up to infinitesimals of arbitrary high order. For that, we identify 
an ideal of trace-class
operators containing all commutation defects; these defects vanish in
the classical case. An appropriately modified first-order condition is
given, which distinguishes Dirac operators with classical spectra.

A discussion of the Connes--Moscovici local index formula for the spectral geometry presented in this paper is currently under investigation 
and will be soon reported elsewhere.

\section{Algebraic preliminaries}
\label{sec:alg-defns}

\begin{defn}
\label{df:SUq-two}
Let $q$ be a real number with $0 < q < 1$, and let
$\A = \A(SU_q(2))$ be the $*$-algebra generated by $a$ and~$b$,
subject to the following commutation rules:
\begin{gather}
ba = q ab,  \qquad  b^*a = qab^*, \qquad bb^* = b^*b,
\nn \\
a^*a + q^2 b^*b = 1,  \qquad  aa^* + bb^* = 1.
\label{eq:suq2-relns}
\end{gather}
As a consequence, $a^*b = q ba^*$ and $a^*b^* = q b^* a^*$. This
becomes a Hopf $*$-algebra under the coproduct
\begin{align*}
\cop a &:= a \ox a - q\,b \ox b^*,
\\
\cop b &:= b \ox a^* + a \ox b,
\end{align*}
counit $\eps(a) = 1$, $\eps(b) = 0$, and antipode $Sa = a^*$,
$Sb = - qb$, $Sb^* = - q^{-1}b^*$, $Sa^* = a$.
\end{defn}

\begin{rem}
Here we follow Majid's ``lexicographic convention''
\cite{MajidPinkBook,MajidStwo} (where, with $c = -qb^*$, $d = a^*$, a
factor of $q$ is needed to restore alphabetical order). Another
much-used convention is related to ours by
$a \otto a^*$, $b \otto -b$; see, for instance, 
\cite{ChakrabortyPEqvt,ConnesSUq}.
\end{rem}

\begin{defn}
\label{df:Uqsu2}
The Hopf $*$-algebra $\U = \U_q(su(2))$ is generated as an algebra
by elements $e,f,k$, with $k$ invertible, satisfying the relations
\begin{equation}
ek = qke,  \qquad  kf = qfk, \qquad
k^2 - k^{-2} = (q - q^{-1})(fe - ef),
\label{eq:Uqsu2-relns}
\end{equation}
and its coproduct $\cop$ is given by
$$
\cop k = k \ox k,  \qquad
\cop e = e \ox k + k^{-1} \ox e,  \qquad
\cop f = f \ox k + k^{-1} \ox f.
$$
Its counit $\epsilon$, antipode $S$, and star structure $^*$ are given
respectively by
\begin{align*}
\epsilon(k) &= 1,  &  Sk &= k^{-1},     &  k^* &= k,
\nn \\
\epsilon(f) &= 0,  &  Sf &= -q f,       &  f^* &= e,
\\
\epsilon(e) &= 0,  &  Se &= -q^{-1} e,  &  e^* &= f.
\nn
\end{align*}
There is an automorphism $\vth$ of $\U_q(su(2))$ defined on the
algebra generators by
\begin{equation}
\vth(k) := k^{-1},  \quad  \vth(f) := -e,  \quad  \vth(e) := -f.
\label{eq:U-autom}
\end{equation}
\end{defn}

\begin{rem}
\label{rem:conventions}
We recall that there is another convention for the generators of
$\U_q(su(2))$ in widespread use: see \cite{Kassel}, for instance. The
handy compendium \cite{KlimykS} gives both versions, denoting by
$\breve U_q(su(2))$ the version which we adopt here. However, the
parameter $q$ of this paper corresponds to $q^{-1}$ in~\cite{KlimykS},
or alternatively, we keep the same~$q$ but exchange $e$ and~$f$ of
that book; the equivalence of these procedures is immediate from the
above formulas \eqref{eq:Uqsu2-relns}.

The older literature uses the convention which we follow here, with
generators usually written as $K = k$, $X^+ = f$, $X^- = e$.
\end{rem}

We employ the so-called ``$q$-integers'', defined for each $n \in \Z$
as
\begin{equation}
[n] = [n]_q := \frac{q^n - q^{-n}}{q - q^{-1}}
\sepword{provided}  q \neq 1.
\label{eq:q-integer}
\end{equation}

\begin{defn}
\label{df:pairing}
There is a bilinear pairing between $\U$ and $\A$, defined on
generators by
$$
\<k, a> = q^\half, \quad  \<k, a^*> = q^{-\half}, \quad
\<e, -qb^*> = \<f, b> = 1,
$$
with all other couples of generators pairing to~0. It satisfies
\begin{equation}
\<(Sh)^*, x> = \overline{\<h, x^*>},
\sepword{for all} h \in \U,\ x \in \A.
\label{eq:pair-star}
\end{equation}
{}We regard $\U$ as a subspace of the linear dual of~$\A$ via this
pairing. There are canonical left and right $\U$-module algebra
structures on~$\A$ \cite{WoronowiczCourse} such that
$$
\<g, h \lt x> := \<gh, x>,  \quad  \<g, x \rt h> := \<hg, x>,
\sepword{for all} g,h \in \U,\ x \in \A.
$$
They are given by $h \lt x := (\id \ox h)\,\cop x$ and
$x \rt h := (h \ox \id)\,\cop x$, or equivalently by
\begin{equation}
h \lt x := \co{x}{1} \,\<h, \co{x}{2}>, \qquad
x \rt h := \<h, \co{x}{1}>\, \co{x}{2},
\label{eq:Hacts-lr}
\end{equation}
using the Sweedler notation $\cop x =: \co{x}{1} \ox \co{x}{2}$ with
implicit summation.
\end{defn}

The right and left actions of $\U$ on $\A$ are mutually commuting:
$$
(h \lt a) \rt g = (\co{a}{1} \,\<h, \co{a}{2}>) \rt g
= \<g, \co{a}{1}> \,\co{a}{2}\, \<h, \co{a}{3}>
= h \lt (\<g, \co{a}{1}>\, \co{a}{2}) = h \lt (a \rt g),
$$
and it follows from \eqref{eq:pair-star} that the star structure is
compatible with both actions:
$$
h \lt x^* = ((Sh)^* \lt x)^*,  \quad  x^* \rt h = (x \rt (Sh)^*)^*,
\sepword{for all}  h \in \U, \ x \in \A.
$$

On the generators, the left action is given explicitly by
\begin{align}
k \lt a &= q^\half a,    &  k \lt a^* &= q^{-\half} a^*,  &
k \lt b &= q^{-\half} b, &  k \lt b^* &= q^\half b^*,
\nn \\
f \lt a &= 0,            &  f \lt a^* &= - q b^*,         &
f \lt b &= a,            &  f \lt b^* &= 0,
\label{eq:Uact-left} \\
e \lt a &= b,            &  e \lt a^* &= 0,               &
e \lt b &= 0,            &  e \lt b^* &= - q^{-1} a^*,
\nn
\end{align}
and the right action is likewise given by
\begin{align}
a \rt k &= q^\half a,   &  a^* \rt k &= q^{-\half} a^*,  &
b \rt k &= q^\half b,   &  b^* \rt k &= q^{-\half} b^*,
\nn \\
a \rt f &= -q b^*,      &  a^* \rt f &= 0,               &
b \rt f &= a^*,         &  b^* \rt f &= 0,
\label{eq:Uact-right} \\
a \rt e &= 0,           &  a^* \rt e &= b,               &
b \rt e &= 0,           &  b^* \rt e &= - q^{-1} a.
\nn
\end{align}

We remark in passing that since $\A$ is also a Hopf algebra, the left
and right actions are linked through the antipodes:
$$
S(Sh \lt x) = Sx \rt h.
$$
Indeed, it is immediate from \eqref{eq:Hacts-lr} and the duality
relation $\<Sh, y> = \<h, Sy>$ that
$$
S(Sh \lt x) = S(\co{x}{1})\, \<Sh, \co{x}{2}>
= S(\co{x}{1})\, \<h, S(\co{x}{2})>
= \co{(Sx)}{2}\, \<h, \co{(Sx)}{1}> = Sx \rt h.
$$

As noted in \cite{GoverZ}, for instance, the invertible antipode
of~$\U$ serves to transform the right action $\rt$ into a second left
action of~$\U$ on~$\A$, commuting with the first. Here we also use the
automorphism $\vth$ of~\eqref{eq:U-autom}, and define
$$
h \. x := x \rt S^{-1}(\vth(h)).
$$
Indeed, it is immediate that
$$
g \. (h \. x) = (x \rt S^{-1}(\vth h)) \rt S^{-1}(\vth g)
= x \rt (S^{-1}(\vth h) S^{-1}(\vth g)) = x \rt (S^{-1}(\vth(gh))
= gh \. x,
$$
i.e., it is a left action. We tabulate this action directly from
\eqref{eq:Uact-right}:
\begin{align}
k \. a &= q^\half a,   &  k \. a^* &= q^{-\half} a^*,  &
k \. b &= q^\half b,   &  k \. b^* &= q^{-\half} b^*,
\nn \\
f \. a &= 0,           &  f \. a^* &= q b,             &
f \. b &= 0,           &  f \. b^* &= - a,
\label{eq:Uact-newl} \\
e \. a &= - b^*,       &  e \. a^* &= 0,               &
e \. b &= q^{-1} a^*,  &  e \. b^* &= 0.
\nn
\end{align}

In the ``classical'' case $q = 1$, we use the well-known 
identifications $SU(2) \approx \Sf^3 \approx \Spin(4)/\Spin(3)
= (SU(2) \x SU(2))/SU(2)$; on quotienting out the diagonal $SU(2)$
subgroup of $\Spin(4)$, we realize $SU(2)$ as the base space of the 
principal spin bundle $\Spin(4) \to \Sf^3$, with projection map
$(g,h) \mapsto gh^{-1}$. The action of $\Spin(4)$ on $SU(2)$ is given
by $(g,h) \. x := gxh^{-1}$, and the stabilizer of~$1$ is the diagonal
$SU(2)$ subgroup. We may choose to regard this as a pair of commuting
actions of $SU(2)$ on the base space $SU(2)$, apart from the nuance of
switching one of them from a right to a left action via the group
inversion map. The foregoing pair of actions of~$\U_q(su(2))$
on~$\A(SU_q(2))$ extends this scheme to the case $q \neq 1$.

\vspace{6pt}

We recall \cite{KlimykS} that $\A$ has a vector-space basis
consisting of matrix elements of its irreducible corepresentations,
$\set{t^l_{mn} : 2l \in \N,\ m,n = -l,\dots,l-1,l}$, where
$$ 
t^0_{00} = 1, \qquad
t^\half_{\half,\half} = a, \qquad
t^\half_{\half,-\half} = b.
$$ 
The coproduct has the matricial form
$\cop t^l_{mn} = \sum_k t^l_{mk} \ox t^l_{kn}$, while the product is
given by
\begin{equation}
t^j_{rs} t^l_{mn} = \sum_{k=|j-l|}^{j+l}
\CGq jlk rm{r+m} \CGq jlk sn{s+n} t^k_{r+m,s+n},
\label{eq:matelt-prod}
\end{equation}
where the $C_q(-)$ factors are $q$-Clebsch--Gordan coefficients
\cite{BiedenharnLo,KirillovR}.

The Haar state on the $C^*$-completion $C(SU_q(2))$, which we shall
denote by~$\psi$, is faithful, and it is determined by setting
$\psi(1) := 1$ and $\psi(t^l_{mn}) := 0$ if $l > 0$. (The Haar state
is usually denoted by $h$, but here we use $h$ for a generic element
of $\U$ instead.) Let $\H_\psi = L^2(SU_q(2),\psi)$ be the Hilbert
space of its GNS representation; then the GNS map
$\eta\: C(SU_q(2)) \to \H_\psi$ is injective and satisfies
\begin{equation}
\|\eta(t^l_{mn})\|^2 = \psi((t^l_{mn})^* \, t^l_{mn})
= \frac{q^{-2m}}{[2l+1]},
\label{eq:matelt-norm}
\end{equation}
and the vectors $\eta(t^l_{mn})$ are mutually orthogonal. From the
formula
$$
\CGq ll0 {-m}m0 = (-1)^{l+m} \frac{q^{-m}}{[2l+1]^\half},
$$
we see that the involution in $C(SU_q(2))$ is given by
\begin{equation}
(t^l_{mn})^* = (-1)^{2l+m+n} q^{n-m} \, t^l_{-m,-n}.
\label{eq:matelt-star}
\end{equation}
In particular, $t^\half_{-\half,\half} = -q b^*$ and 
$t^\half_{-\half,-\half} = a^*$, as expected.

An orthonormal basis of $\H_\psi$ is obtained by normalizing the
matrix elements, using~\eqref{eq:matelt-norm}:
\begin{equation}
\ket{lmn} := q^m \,[2l+1]^\half \,\eta(t^l_{mn}).
\label{eq:matelt-onb}
\end{equation}

\section{Equivariant representation of {\boldmath$\A(SU_q(2))$}}
\label{sec:eqvt-repns}

Let $\U$ be a Hopf algebra and let $\A$ be a left $\U$-module algebra.
A representation of $\A$ on a vector space $V$ is called
$\U$-equivariant if there is also an algebra representation of $\U$
on~$V$, satisfying the following compatibility relation:
$$
h(x\xi) = (\co{h}{1} \lt x)(\co{h}{2} \xi),  \qquad
h \in \U,\ x \in \A,\ \xi \in V,
$$
where $\lt$ denotes the Hopf action of $\U$ on~$\A$. If $\A$ is
instead a right $\U$-module algebra, the appropriate compatibility
relation is $x(h\xi) = \co{h}{1}((x \rt \co{h}{2})\xi)$. Also, if $\A$
is an $\U$-bimodule algebra (carrying commuting left and right Hopf
actions of~$\U$), one can demand both of these conditions
simultaneously for pair of representations of $\A$ and~$\U$ on the
same vector space~$V$.

In the present case, it turns out to be simpler to consider
equivariance under two commuting left Hopf actions, as exemplified in
the previous section. We shall first work out in detail a construction
of the regular representation of the Hopf algebra $\A(SU_q(2))$,
showing how it is determined by its equivariance properties.

\vspace{6pt}

We begin with the known representation theory~\cite{KlimykS} of
$\U_q(su(2))$. The irreducible finite dimensional representations
$\sg_l$ of $\U_q(su(2))$ are labelled by nonnegative half-integers
$l = 0,\half,1,\sesq,2,\dots$, and they are given by
\begin{align}
\sg_l(k)\,\ket{lm} &= q^m \,\ket{lm},
\nn \\
\sg_l(f)\,\ket{lm} &= \sqrt{[l-m][l+m+1]} \,\ket{l,m+1},
\label{eq:uqsu2-repns} \\
\sg_l(e)\,\ket{lm} &= \sqrt{[l-m+1][l+m]} \,\ket{l,m-1},
\nn
\end{align}
where the vectors $\ket{lm}$, for $m = -l, -l+1,\dots, l-1, l$, form a
basis for the irreducible $\U$-module $V_l$, and the brackets denote
$q$-integers as in~\eqref{eq:q-integer}. Moreover, $\sg_l$ is a
$*$-representation of $\U_q(su(2))$, with respect to the hermitian
scalar product on $V_l$ for which the vectors $\ket{lm}$ are
orthonormal.

\begin{rem}
\label{rem:uqsu2-repns}
The irreducible representations \eqref{eq:uqsu2-repns} coincide with
those of $\breve U_q(su(2))$ in \cite{KlimykS}, after exchange of $e$
and~$f$ (see Remark~\ref{rem:conventions}). Further results on the
representation theory of $\U_q(su(2))$ are taken from
\cite[Chap.~3]{KlimykS} without comment; in particular we use the
$q$-Clebsch--Gordan coefficients found therein for the decomposition
of tensor product representations. An alternative source for these
coefficients is \cite{BiedenharnLo}, although their $q^\half$ is
our~$q$.
\end{rem}

\begin{defn}
\label{df:Heqvt-repn}
Let $\la$ and $\rho$ be mutually commuting representations of the Hopf
algebra $\U$ on a vector space~$V$. A representation $\pi$ of the
$*$-algebra $\A$ on $V$ is \textit{$(\la,\rho)$-equivariant} if the
following compatibility relations hold:
\begin{align}
\la(h)\,\pi(x) \xi &= \pi(\co{h}{1} \. x) \, \la(\co{h}{2}) \xi,
\nn \\
\rho(h)\,\pi(x)\xi &= \pi(\co{h}{1} \lt x)\,\rho(\co{h}{2}) \xi,
\label{eq:covar-repn}
\end{align}
for all $h \in \U$, $x \in \A$ and $\xi \in V$.
\end{defn}

\vspace{6pt}

We shall now exhibit an equivariant representation of $\A(SU_q(2))$ on
the preHilbert space which is the (algebraic) direct sum
$$
V := \bigoplus_{2l=0}^\infty  V_l \ox V_l.
$$
The two $\U_q(su(2))$ symmetries $\la$ and $\rho$ will act on the
first and the second leg of the tensor product respectively; both
actions will be via the irreps \eqref{eq:uqsu2-repns}. In other words,
$$
\la(h) = \sg_l(h) \ox \id,  \qquad  \rho(h) = \id \ox \sg_l(h)
\qquad\text{on } V_l \ox V_l.
$$
We abbreviate $\ket{lmn} := \ket{lm} \ox \ket{ln}$, for
$m,n = -l,\dots,l-1,l$; these form an orthonormal basis for
$V_l \ox V_l$, for each fixed~$l$. (As we shall see, this is
consistent with our labelling \eqref{eq:matelt-onb} of the orthonormal
basis of $\H_\psi$ in the previous section.) Also, we adopt a
shorthand notation:
$$
l^\pm := l \pm \half,  \quad
m^\pm := m \pm \half,  \quad
n^\pm := n \pm \half.
$$

\begin{prop}
\label{pr:suq2-repn}
A $(\la,\rho)$-equivariant $*$-representation $\pi$ of $\A(SU_q(2))$
on the Hilbert space $V$ of~\eqref{pr:suq2-repn} must have the
following form:
\begin{align}
\pi(a) \,\ket{lmn}
&= A^+_{lmn} \ket{l^+ m^+ n^+} + A^-_{lmn} \ket{l^- m^+ n^+},
\nn \\
\pi(b) \,\ket{lmn}
&= B^+_{lmn} \ket{l^+ m^+ n^-} + B^-_{lmn} \ket{l^- m^+ n^-},
\nn \\
\pi(a^*) \,\ket{lmn}
&= \Onda A^+_{lmn} \ket{l^+m^-n^-} + \Onda A^-_{lmn} \ket{l^-m^-n^-},
\label{eq:reg-repn} \\
\pi(b^*) \,\ket{lmn}
&= \Onda B^+_{lmn} \ket{l^+m^-n^+} + \Onda B^-_{lmn} \ket{l^-m^-n^+},
\nn
\end{align}
where the constants $A^\pm_{lmn}$ and $B^\pm_{lmn}$ are, up to phase
factors depending only on~$l$, given by
\begin{align}
A^+_{lmn} &= q^{(-2l+m+n-1)/2}
\biggl( \frac{[l+m+1][l+n+1]}{[2l+1][2l+2]} \biggr)^\half \!,
\nn \\
A^-_{lmn} &= q^{(2l+m+n+1)/2}
\biggl( \frac{[l-m][l-n]}{[2l][2l+1]} \biggr)^\half \!,
\nn \\
B^+_{lmn} &= q^{(m+n-1)/2}
\biggl( \frac{[l+m+1][l-n+1]}{[2l+1][2l+2]} \biggr)^\half \!,
\label{eq:reg-repn-coeffs} \\
B^-_{lmn} &= - q^{(m+n-1)/2}
\biggl( \frac{[l-m][l+n]}{[2l][2l+1]} \biggr)^\half \!,
\nn
\end{align}
and the other coefficients are complex conjugates of these, namely,
\begin{equation}
\Onda A^\pm_{lmn} = (A^\mp_{l^\pm m^- n^-})^\star,  \qquad
\Onda B^\pm_{lmn} = (B^\mp_{l^\pm m^- n^+})^\star.
\label{eq:conj-coeff}
\end{equation}
\end{prop}

\begin{proof}
First of all, notice that hermiticity of $\pi$ entails the relations
\eqref{eq:conj-coeff}. We now use the covariance properties
\eqref{eq:covar-repn}. When $h = k$, they simplify to
\begin{equation}
\la(k)\, \pi(x)\,\xi = \pi(k \. x)\, \la(k)\,\xi,  \qquad
\rho(k)\,\pi(x)\,\xi = \pi(k \lt x)\,\rho(k)\,\xi.
\label{eq:covar-k}
\end{equation}
Thus, for instance, when $x = a$ we find the relations
\begin{align*}
\la(k)\,\pi(a)\,\ket{lmn}
&= \pi(q^\half a) \bigl( q^m\ket{lmn} \bigr)
= q^{m+\half} \pi(a)\,\ket{lmn},
\\
\rho(k)\,\pi(a)\,\ket{lmn}
&= \pi(q^{\half}a) \bigl( q^n\ket{lmn} \bigr)
= q^{n+\half} \pi(a)\,\ket{lmn},
\end{align*}
where we have invoked $k \. a = k \lt a = q^\half a$. We conclude that
$\pi(a)\,\ket{lmn}$ must lie in the closed span of the basis vectors
$\ket{l' m^+ n^+}$. A similar argument with $x = b$
in~\eqref{eq:covar-k} shows that $\pi(b)$ increments $n$ and
decrements $m$ by~$\half$, since $k \. b = q^\half b$ while
$k \lt b = q^{-\half} b$. The analogous behaviour for $x = a^*$ and
$x = b^*$ follows in the same way from \eqref{eq:Uact-left} and
\eqref{eq:Uact-newl}.

Thus, $\pi(a)\ket{lmn}$ is a (possibly infinite) sum
\begin{equation}
\pi(a)\,\ket{lmn} = \tsum_{l'} C_{l'lmn}\, \ket{l' m^+ n^+},
\label{eq:pi-a-xpan}
\end{equation}
where the sum runs over nonnegative half-integers
$l' = 0, \half, 1, \sesq, \dots$.

Next, we call on \eqref{eq:covar-repn} with $h = f$, $x = a$, to get
$$
\la(f)\,\pi(a) \xi
= \pi(f \. a)\,\la(k)\xi + \pi(k^{-1} \. a)\,\la(f)\xi
= q^{-\half} \pi(a)\,\la(f)\xi,
$$
on account of \eqref{eq:Uact-left}. Consequently,
$\la(f)^r \pi(a) = q^{-r/2} \pi(a)\,\la(f)^r$ for $r = 1,2,3,\dots$.
On applying $\la(f)^r$ to both sides of \eqref{eq:pi-a-xpan}, we
obtain on the left hand side a multiple of $\pi(a)\,\ket{l,m+r,n}$,
which vanishes for $m + r > l$; and on the right hand side we get
$\tsum_{l'} C_{l'lmn}\,D_{l'mr}\, \ket{l', m^+ + r, n^+}$, where
$D_{l'mr} \neq 0$ as long as $m + r + \half \leq l'$. We conclude that
$C_{l'lmn} = 0$ for $l' > l + \half$, by linear independence of these
summands.

To get a lower bound on the range of the index $l'$ in
\eqref{eq:pi-a-xpan}, we consider the analogous expansion
$\pi(a^*)\,\ket{lmn} = \tsum_{l'} \Onda C_{l'lmn}\, \ket{l' m^- n^-}$.
Now
$\la(e)^r \pi(a^*)\,\ket{lmn} = q^{r/2} \pi(a^*)\,\la(e)^r\,\ket{lmn} 
\allowbreak \propto \pi(a^*)\,\ket{l,m-r,n}$ vanishes for
$m - r < -l$; while
$\la(e)^r \ket{l' m^- n^-} = F_{l'mr}\, \ket{l', m^- - r, n^-}$ with
$F_{l'mr} \neq 0$ for $m - r - \half \geq -l'$. Again we conclude that
$\Onda C_{l'lmn} = 0$ for $l' > l + \half$. However, since $\pi$ is a
$*$-representation, the matrix element $\braCket{l'm'n'}{\pi(a)}{lmn}$
is the complex conjugate of $\braCket{lmn}{\pi(a^*)}{l'm'n'}$, which
vanishes for $l > l' + \half$, so that the indices in
\eqref{eq:pi-a-xpan} satisfy $l - \half \leq l' \leq l + \half$.
Clearly, $l' = l$ is ruled out because $l - m$ and $l' - m \pm \half$
must both be integers.

Therefore, $\pi(a)$ and also $\pi(a^*)$ have the structure indicated
in \eqref{eq:reg-repn}. A parallel argument shows the corresponding
result for $\pi(b)$ and~$\pi(b^*)$.

The coefficients which appear in \eqref{eq:reg-repn-coeffs} may be
determined by further application of the equivariance relations.
Since $f \lt a = 0$ and $e \lt b = 0$, then by applying $\rho(f)$ and
$\rho(e)$ to the first two relations of \eqref{eq:reg-repn}, we obtain
the following recursion relations for the coefficients $A^\pm_{lmn}$,
$B^\pm_{lmn}$:
\begin{align*}
A^+_{lmn} [l+n+2]^\half &= q^{-\half} A^+_{lm,n+1} [l+n+1]^\half,
\\
A^-_{lmn} [l-n-1]^\half &= q^{-\half} A^-_{lm,n+1} [l-n]^\half, 
\\
B^+_{lmn} [l-n+2]^\half &= q^\half B^+_{lm,n-1} [l-n+1]^\half, 
\\
B^-_{lmn} [l+n-1]^\half &= q^\half B^-_{lm,n-1} [l+n]^\half.
\end{align*}
Then, applying $\la(f)$ to the same pair of equations, we further find
that
\begin{subequations}
\label{eq:l-only}
\begin{align}
A^+_{lmn} [l+m+2]^\half &= q^{-\half} A^+_{l,m+1,n} [l+m+1]^\half,
\nn \\
A^-_{lmn} [l-m-1]^\half &= q^{-\half} A^-_{l,m+1,n} [l-m]^\half, 
\nn \\
B^+_{lmn} [l+m+2]^\half &= q^{-\half} B^+_{l,m+1,n} [l+m+1]^\half,
\label{eq:lonly-recurs} \\
B^-_{lmn} [l-m-1]^\half &= q^{-\half} B^-_{l,m+1,n} [l-m]^\half.
\nn
\end{align}
These recursions are explicitly solved by
\begin{align}
A^+_{lmn} &= q^{(m+n)/2} [l+m+1]^\half [l+n+1]^\half \ a^+_l, 
\nn \\
A^-_{lmn} &= q^{(m+n)/2} [l-m]^\half [l-n]^\half     \ a^-_l, 
\nn \\
B^+_{lmn} &= q^{(m+n)/2} [l+m+1]^\half [l-n+1]^\half \ b^+_l, 
\label{eq:lonly-coeff} \\
B^-_{lmn} &= q^{(m+n)/2} [l-m]^\half [l+n]^\half     \ b^-_l,
\nn
\end{align}
\end{subequations}
where $a^\pm_l$, $b^\pm_l$ depend only on~$l$.

Once more, we apply the equivariance relations \eqref{eq:covar-repn};
this time, we use
\begin{equation}
\rho(e) \pi(a) = \pi(e \lt a) \rho(k) + \pi(k^{-1} \lt a) \rho(e)
= \pi(b) \rho(k) + q^{-\half} \pi(a) \rho(e).
\label{eq:ab-link}
\end{equation}
Applied to $\ket{lmn}$, it yields an equation between linear
combinations of $\ket{l^+ m^+ n^-}$ and $\ket{l^- m^+ n^-}$; equating
coefficients, we find
$$
b^+_l = q^l a^+_l,  \qquad  b^-_l = - q^{-l-1} a^-_l.
$$
Furthermore, applying also to $\ket{lmn}$ the relation
\begin{align}
\la(e) \pi(b) &= \pi(e \. b) \la(k) + \pi(k^{-1} \. b) \la(e) 
\nn \\
&= q^{-1} \pi(a^*) \la(k) + q^{-\half} \pi(b) \la(e),
\label{eq:bac-link}
\end{align}
we get, after a little simplification and use of 
\eqref{eq:conj-coeff},
$$
(a^-_{l+\half})^\star = q^{2l+\sesq} \, a^+_l.
$$

It remains only to determine the parameters $a^+_l$. We turn to the
algebra commutation relation $ba = qab$ and compare coefficients in
the expansion of
$\pi(b)\pi(a)\,\ket{lmn} = q\,\pi(a)\pi(b)\,\ket{lmn}$. Those of
$\ket{l+1,m+1,n}$ and $\ket{l-1,m+1,n}$ already coincide; but from the
$\ket{l,m+1,n}$ terms, we get the identity
$$
q [2l+2]\, |a^+_l|^2 = [2l]\, |a^+_{l-\half}|^2.
$$
This can be solved immediately, to give
$$
a^+_l = \frac{C \zeta_l\,q^{-l}}{[2l+1]^\half[2l+2]^\half},
$$
where $C$ is a positive constant, and $\zeta_l$ is a phase factor
which can be absorbed in the basis vectors $\ket{lmn}$; hereinafter we
take $\zeta_l = 1$ (we comment on that choice at the end of the
section).

Finally, from the relation $a^*a + q^2 b^*b = 1$ we obtain
$$
1 = \braCket{000}{\pi(a^*a + q^2 b^*b)}{000} 
= |a^+_0|^2 + q^2 |b^+_0|^2 = (1 + q^2) C^2 /[2] = q\,C^2,
$$
and thus $C = q^{-\half}$. We therefore find that
\begin{align*}
a^+_l &= \frac{q^{-l-\half}}{[2l+1]^\half[2l+2]^\half}, &
a^-_l &= \frac{q^{l+\half}}{[2l]^\half[2l+1]^\half},    \nn \\
b^+_l &= \frac{q^{-\half}}{[2l+1]^\half[2l+2]^\half},   &
b^-_l &= - \frac{q^{-\half}}{[2l]^\half[2l+1]^\half},
\end{align*}
and substitution in \eqref{eq:lonly-coeff} yields the coefficients
\eqref{eq:reg-repn-coeffs}.
\end{proof}

It is easy to check that the formulas \eqref{eq:reg-repn} give
precisely the left regular representation $\pi_\psi$ of $\A(SU_q(2))$.
Indeed, that representation was implicitly given already by the
product rule \eqref{eq:matelt-prod}. From \cite[(3.53)]{BiedenharnLo}
we obtain
\begin{align}
\CGq {\half}{l}{l^+} {\half}{m}{m^+}
&= q^{-\half(l-m)} \,\frac{[l+m+1]^\half}{[2l+1]^\half}, 
\nn \\
\CGq {\half}{l}{l^+} {-\half}{m}{m^-}
&= q^{\half(l+m)} \,\frac{[l-m+1]^\half}{[2l+1]^\half},
\nn \\
\CGq {\half}{l}{l^-} {\half}{m}{m^+}
&= q^{\half(l+m+1)} \,\frac{[l-m]^\half}{[2l+1]^\half}, 
\label{eq:CG-coeffs} \\
\CGq {\half}{l}{l^-} {-\half}{m}{m^-}
&= - q^{-\half(l-m+1)} \,\frac{[l+m]^\half}{[2l+1]^\half}.
\nn
\end{align}
By setting $j = r = s = \half$ in \eqref{eq:matelt-prod}, we find
$$
\pi_\psi(a) \eta(t^l_{mn})
= \sum_\pm \CGq {\half}{l}{l^\pm} {\half}{m}{m^+}
\CGq {\half}{l}{l^\pm} {\half}{n}{n^+} \eta(t^{l^\pm}_{m^+n^+}).
$$
Taking the normalization \eqref{eq:matelt-onb} into account, this
becomes
\begin{align*}
\pi_\psi(a) \ket{lmn}
&= q^{-\half} \frac{[2l+1]^\half}{[2l+2]^\half}
\CGq {\half}{l}{l^+} {\half}{m}{m^+} 
\CGq {\half}{l}{l^+} {\half}{n}{n^+} \,\ket{l^+m^+n^+}
\\
&\qquad + q^{-\half} \frac{[2l+1]^\half}{[2l]^\half}
\CGq {\half}{l}{l^-} {\half}{m}{m^+} 
\CGq {\half}{l}{l^-} {\half}{n}{n^+} \,\ket{l^-m^+n^+}
\\
&= q^{\half(-2l+m+n-1)} \,
\frac{[l+m+1]^\half[l+n+1]^\half}{[2l+1]^\half[2l+2]^\half} \,
\ket{l^+m^+n^+}
\\
&\qquad + q^{\half(2l+m+n+1)} \,
\frac{[l-m]^\half[l-n]^\half}{[2l]^\half[2l+1]^\half}\,\ket{l^-m^+n^+}
\\
&= \pi(a) \ket{lmn}.
\end{align*}
A similar calculation, using \eqref{eq:CG-coeffs} again, shows that 
$\pi(b) = \pi_\psi(b)$. Since $a$ and $b$ generate $\A$ as a 
$*$-algebra, we conclude that $\pi = \pi_\psi$. (It should be noted 
that $\pi_\psi$ has already been exhibited in \cite{ChakrabortyPEqvt}
in the same way, albeit with different convention for the algebra 
generators.)

The identification~\eqref{eq:matelt-onb} embeds the prehilbert space
$V$ densely in the Hilbert space $\H_\psi$, and the representation
$\pi_\psi$ extends to the GNS representation of $C(SU_q(2))$
on~$\H_\psi$, as described by the Peter-Weyl theorem
\cite{KlimykS,WoronowiczCourse}. In like manner, all other 
representations of~$\A$ exhibited in this paper extend to 
$C^*$-algebra representations of $C(SU_q(2))$ on the appropriate 
Hilbert spaces.

\vspace{6pt}

The only lack of uniqueness in the proof of
Proposition~\ref{pr:suq2-repn} involved the choice of
the phase factors~$\zeta_l$; if $Z$ is the linear operator on~$V$
which multiplies vectors in $V_l \ox V_l$ by~$\zeta_l$, then $Z$
commutes with each $\la(h)$ and $\rho(g)$, and extends to a unitary
operator on $\H_\psi$. In other words, any $(\la,\rho)$-equivariant
representation $\pi$ extends to~$\H_\psi$ and is unitarily equivalent
to the left regular representation. The (standard) choice
$\zeta_l = 1$ ensures that \emph{all coefficients $A^\pm_{lmn}$ and
$B^\pm_{lmn}$ are real}: it is indeed an extension of the
Conden-Shortley phase convention~\cite{BiedenharnLk}.

\section{The spin representation}
\label{sec:spinor-repn}

The left regular representation $\pi$ of $\A$, constructed in the
previous section, can be amplified to $\pi' = \pi \ox \id$ on
$V \ox \C^2$. In the commutative case when $q = 1$, this yields the
spinor representation of $SU(2)$, because the spinor bundle is
parallelizable: $S \simeq SU(2) \x \C^2$, although one needs to
specify the trivialization. The representation theory of $\U$ (and the
corepresentation theory of~$\A$) follows the same pattern; only the
Clebsch--Gordan coefficients need to be modified \cite{KirillovR} when
$q \neq 1$.

To fix notations, we take
$$
W := V \ox \C^2 = V \ox V_\half,
$$
and its Clebsch--Gordan decomposition is the (algebraic) direct sum
\begin{equation}
W = \biggl( \bigoplus_{2l=0}^\infty V_l \ox V_l \biggr) \ox V_\half
\simeq V_\half \oplus \bigoplus_{2j=1}^\infty
(V_{j+\half} \ox V_j) \oplus (V_{j-\half} \ox V_j).
\label{eq:CG-decomp}
\end{equation}
We rename the finite-dimensional spaces on the right hand side as
\begin{equation}
W = W_0^\up \oplus \bigoplus_{2j\geq 1} W_j^\up \oplus W_j^\dn ,
\label{eq:spinor-decomp}
\end{equation}
where $W_j^\up \simeq V_{j+\half} \ox V_j$ and
$W_j^\dn \simeq V_{j-\half} \ox V_j$, so that 
\begin{equation}
\begin{aligned}
\dim W_j^\up &= (2j + 1)(2j + 2), 
\\
\dim W_j^\dn &= 2j(2j + 1),
\end{aligned}
\quad
\begin{aligned}
\text{for } j &= 0,\half,1,\sesq,\dots~,
\\
\text{for } j &= \half,1,\sesq,\dots~.
\end{aligned}
\label{eq:Wj-dims}
\end{equation}

\begin{defn}
\label{df:Heqvt-spinor}
We amplify the representation $\rho$ of $\U$ on~$V$ to 
$\rho' = \rho \ox \id$ on $W = V \ox \C^2$. However, we replace 
$\la$ on~$V$ by its tensor product with $\sg_\half$ on~$\C^2$:
$$
\la'(h) := (\la \ox \sg_\half)(\cop h)
= \la(\co{h}{1}) \ox \sg_\half(\co{h}{2}).
$$
It is straightforward to check that the representations $\la'$ and 
$\rho'$ on~$W$ commute, and that the representation $\pi'$ of $\A$ on 
$W$ is $(\la',\rho')$-equivariant:
\begin{align}
\la'(h)\,\pi'(x) \psi &= \pi'(\co{h}{1} \. x) \, \la'(\co{h}{2}) \psi,
\nn \\
\rho'(h)\,\pi'(x)\psi &= \pi'(\co{h}{1} \lt x)\,\rho'(\co{h}{2}) \psi,
\label{eq:covar-spin}
\end{align}
for all $h \in \U$, $x \in \A$ and $\psi \in W$.
\end{defn}

To determine an explicit basis for $W$ which is well-adapted to 
$(\la',\rho')$-equivariance, consider the following vectors in
$V \ox \C^2$:
\begin{align*}
\phantom{-} c_{lm} \,\ket{lmn} \ox \ket{\half,-\half}
&+ s_{lm} \,\ket{l,m-1,n} \ox \ket{\half,+\half},
\\[\jot]
- s_{lm} \,\ket{lmn} \ox \ket{\half,-\half}
&+ c_{lm} \,\ket{l,m-1,n} \ox \ket{\half,+\half},
\end{align*}
where 
$$
c_{lm} := q^{-(l+m)/2}  \,\frac{[l-m+1]^\half}{[2l+1]^\half}, \qquad
s_{lm} := q^{(l-m+1)/2} \,\frac{[l+m]^\half}{[2l+1]^\half}
$$
are the $q$-Clebsch--Gordan coefficients corresponding to the above
decomposition \eqref{eq:CG-decomp}, satisfying
$c_{lm}^2 + s_{lm}^2 = 1$. These are eigenvectors for $\la'(C_q)$,
where $C_q := q k^2 + q^{-1} k^{-2} + (q - q^{-1})^2 ef$ is the
Casimir element of $\U$, with respective eigenvalues
$q^{2l+2} + q^{-2l-2}$ and $q^{2l} + q^{-2l}$. Thus, to get a good
basis, one should offset the index~$l$ by $\pm\half$ (as is also
suggested by the decomposition \eqref{eq:spinor-decomp} of~$W$).

For $j = l + \half$, $\mu = m - \half$, with $\mu = -j,\dots,j$ and
$n = -j^-,\dots,j^-$, let
\begin{subequations}
\label{eq:spinor-basis}
\begin{align}
\ket{j\mu n\dn}
&:= C_{j\mu} \,\ket{j^- \mu^+ n} \ox \ket{\half,-\half}
+ S_{j\mu} \,\ket{j^- \mu^- n} \ox \ket{\half,+\half};
\label{eq:spinor-basis-dn}
\\
\intertext{and for $j = l - \half$, $\mu = m - \half$, with 
$\mu = -j,\dots,j$ and $n = -j^+,\dots,j^+$, let}
\ket{j\mu n\up}
&:= - S_{j+1,\mu} \,\ket{j^+ \mu^+ n} \ox \ket{\half,-\half}
+ C_{j+1,\mu} \,\ket{j^+ \mu^- n} \ox \ket{\half,+\half},
\label{eq:spinor-basis-up}
\end{align}
where the coefficients are now
\begin{equation}
C_{j\mu} := q^{-(j+\mu)/2}\,\frac{[j-\mu]^\half}{[2j]^\half},  \qquad
S_{j\mu} := q^{(j-\mu)/2} \,\frac{[j+\mu]^\half}{[2j]^\half}.
\label{eq:spinor-basis-coeffs}
\end{equation}
\end{subequations}
Notice that there are no $\dn$ vectors for $j = 0$. It is now
straightforward, though tedious, to verify that these vectors are
orthonormal bases for the respective subspaces $W_j^\dn$
and~$W_j^\up$.

The Hilbert space of spinors is $\H := \H_\psi \ox \C^2$, which is
just the completion of the algebraic direct sum
\eqref{eq:spinor-decomp}. We may decompose it as
$\H = \H^\up \oplus \H^\dn$, where $\H^\up$ and $\H^\dn$ are the 
respective completions of $\bigoplus_{2j\geq 0} W_j^\up$ and 
$\bigoplus_{2j\geq 1} W_j^\dn$.

\begin{lem}
\label{lm:eigen-down}
The basis vectors $\ket{j\mu n\up}$ and $\ket{j\mu n\dn}$ are joint
eigenvectors for $\la'(k)$ and $\rho'(k)$, and $e,f$ are represented
on them as ladder operators:
\begin{subequations}
\label{eq:spinor-ladder}
\begin{equation}
\begin{aligned}
\la'(k)\ket{j\mu n\up} &= q^\mu \ket{j\mu n\up}, \\
\la'(k)\ket{j\mu n\dn} &= q^\mu \ket{j\mu n\dn},
\end{aligned} \qquad
\begin{aligned}
\rho'(k)\ket{j\mu n\up} &= q^n \ket{j\mu n\up},  \\
\rho'(k)\ket{j\mu n\dn} &= q^n \ket{j\mu n\dn}.
\end{aligned}
\label{eq:spinor-ladder-k}
\end{equation}
Moreover,
\begin{align}
\la'(f)\ket{j\mu n\up}
&= [j - \mu]^\half [j + \mu + 1]^\half \ket{j,\mu+1,n\up},
\nn \\[\jot]
\la'(e)\ket{j\mu n\up}
&= [j + \mu]^\half [j - \mu + 1]^\half \ket{j,\mu-1,n\up},
\nn \\[\jot]
\la'(f)\ket{j\mu n\dn}
&= [j - \mu]^\half [j + \mu + 1]^\half \ket{j,\mu+1,n\dn},
\label{eq:spinor-ladder-L} \\[\jot]
\la'(e)\ket{j\mu n\dn}
&= [j + \mu]^\half [j - \mu + 1]^\half \ket{j,\mu-1,n\dn},
\nn
\end{align}
and
\begin{align}
\rho'(f)\ket{j\mu n\up}
&= [j - n + \half]^\half [j + n + \sesq]^\half \ket{j\mu,n+1,\up},
\nn \\[\jot]
\rho'(e)\ket{j\mu n\up}
&= [j + n + \half]^\half [j - n + \sesq]^\half \ket{j\mu,n-1,\up},
\nn \\[\jot]
\rho'(f)\ket{j\mu n\dn}
&= [j - n - \half]^\half [j + n + \half]^\half \ket{j\mu,n+1,\dn},
\label{eq:spinor-ladder-R} \\[\jot]
\rho'(e)\ket{j\mu n\dn}
&= [j + n - \half]^\half [j - n + \half]^\half \ket{j\mu,n-1,\dn}.
\nn
\end{align}
\end{subequations}
\end{lem}

\vspace{6pt}

The representation $\pi'$ can now be computed in the new spinor basis
by conjugating the form of $\pi \ox \id$ found in
Proposition~\ref{pr:suq2-repn} by the basis transformation
\eqref{eq:spinor-basis}. However, it is more instructive to derive
these formulas from the property of $(\la',\rho')$-equivariance.
First, we introduce a handy notation.

\begin{defn}
\label{df:eigen-spinors}
For $j = 0,\half,1,\sesq,\dots$, with $\mu = -j,\dots,j$ and
$n = -j-\half,\dots,j+\half$, we juxtapose the pair of spinors
$$
\kett{j\mu n} := \begin{pmatrix} \ket{j\mu n\up} \\[2\jot]
\ket{j\mu n\dn} \end{pmatrix},
$$
with the convention that the lower component is zero when
$n = \pm(j+\half)$ or $j = 0$. Furthermore, a matrix with scalar
entries,
$$
A = \begin{pmatrix} A_{\up\up} & A_{\up\dn} \\
A_{\dn\up} & A_{\dn\dn} \end{pmatrix},
$$
is understood to act on $\kett{j\mu n}$ by the rule:
\begin{align}
A \ket{j\mu n \up}
&= A_{\up\up} \ket{j\mu n\up} + A_{\dn\up} \ket{j\mu n\dn},
\nn \\
A \ket{j\mu n \dn}
&= A_{\dn\dn} \ket{j\mu n\dn} + A_{\up\dn} \ket{j\mu n\up}.
\label{eq:act-on-arrows}
\end{align}
\end{defn}

\begin{prop}
\label{pr:spin-repn}
The representation $\pi' := \pi \ox \id$ of $\A$ is given by
\begin{align}
\pi'(a) \,\kett{j\mu n}
&= \a^+_{j\mu n} \kett{j^+ \mu^+ n^+}
 + \a^-_{j\mu n} \kett{j^- \mu^+ n^+},
\nn \\[\jot]
\pi'(b) \,\kett{j\mu n}
&= \b^+_{j\mu n} \kett{j^+ \mu^+ n^-}
 + \b^-_{j\mu n} \kett{j^- \mu^+ n^-},
\nn \\[\jot]
\pi'(a^*) \,\kett{j\mu n}
&= \tilde\a^+_{j\mu n} \kett{j^+ \mu^- n^-}
 + \tilde\a^-_{j\mu n} \kett{j^- \mu^- n^-},
\label{eq:spin-repn} \\[\jot]
\pi'(b^*) \,\kett{j\mu n}
&= \tilde\b^+_{j\mu n} \kett{j^+ \mu^- n^+}
 + \tilde\b^-_{j\mu n} \kett{j^- \mu^- n^+},
\nn
\end{align}
where $\a^\pm_{j\mu n}$ and $\b^\pm_{j\mu n}$ are, up to phase factors
depending only on~$j$, the following triangular $2 \x 2$ matrices:

\noindent
\begin{align}
\a^+_{j\mu n} &= q^{(\mu+n-\half)/2} [j + \mu + 1]^\half
\begin{pmatrix}
q^{-j-\half} \, \frac{[j+n+\sesq]^{1/2}}{[2j+2]} & 0 \\[2\jot]
q^\half \,\frac{[j-n+\half]^{1/2}}{[2j+1]\,[2j+2]} &
q^{-j} \, \frac{[j+n+\half]^{1/2}}{[2j+1]}
\end{pmatrix},
\nn \\[2\jot]
\a^-_{j\mu n} &= q^{(\mu+n-\half)/2} [j - \mu]^\half
\begin{pmatrix}
q^{j+1} \, \frac{[j-n+\half]^{1/2}}{[2j+1]} &
- q^\half \,\frac{[j+n+\half]^{1/2}}{[2j]\,[2j+1]} \\[2\jot]
0 & q^{j+\half} \, \frac{[j-n-\half]^{1/2}}{[2j]}
\end{pmatrix},
\nn \\[2\jot]
\b^+_{j\mu n} &= q^{(\mu+n-\half)/2} [j + \mu + 1]^\half
\begin{pmatrix}
\frac{[j-n+\sesq]^{1/2}}{[2j+2]} & 0 \\[2\jot]
- q^{-j-1} \,\frac{[j+n+\half]^{1/2}}{[2j+1]\,[2j+2]} &
q^{-\half} \, \frac{[j-n+\half]^{1/2}}{[2j+1]}
\end{pmatrix},
\label{eq:spin-coeff}
\\[2\jot]
\b^-_{j\mu n} &= q^{(\mu+n-\half)/2} [j - \mu]^\half
\begin{pmatrix}
- q^{-\half} \, \frac{[j+n+\half]^{1/2}}{[2j+1]} &
- q^j \,\frac{[j-n+\half]^{1/2}}{[2j]\,[2j+1]} \\[2\jot]
0 & - \frac{[j+n-\half]^{1/2}}{[2j]}
\end{pmatrix},
\nn
\end{align}
and the remaining matrices are the hermitian conjugates
$$
\tilde\a^\pm_{j\mu n} = (\a^\mp_{j^\pm \mu^- n^-})^\7,  \qquad
\tilde\b^\pm_{j\mu n} = (\b^\mp_{j^\pm \mu^- n^+})^\7.
$$
\end{prop}

\begin{proof}
The proof of Proposition~\ref{pr:suq2-repn} applies with minor
changes. {}From the analogues of \eqref{eq:covar-k} and the relations
$\la'(f) \pi'(a) = q^{-\half} \pi'(a)\,\la'(f)$ and
$\la'(e) \pi'(a^*) = q^\half \pi'(a^*)\,\la'(e)$, applied to the
spinors $\kett{j\mu n}$, together with the formulas
\eqref{eq:spinor-ladder-k} and \eqref{eq:spinor-ladder-L}, we
determine that $\pi'(a)$ has the indicated form, where the
$\a^\pm_{j\mu n}$ are $2 \x 2$ matrices. The other cases of
\eqref{eq:spin-repn} are handled similarly.

To compute these matrices, we again use the commutation relations
of $\la'(f)$ with $\pi'(a)$ and $\pi'(b)$ to establish recurrence 
relations, analogous to \eqref{eq:lonly-recurs}, which yield
\begin{align*}
\a^+_{j\mu n} &= q^{(\mu+n-\half)/2} [j + \mu + 1]^\half\ A^+_{jn}, &
\a^-_{j\mu n} &= q^{(\mu+n-\half)/2} [j - \mu]^\half    \ A^-_{jn}, \\
\b^+_{j\mu n} &= q^{(\mu+n-\half)/2} [j + \mu + 1]^\half\ B^+_{jn}, &
\b^-_{j\mu n} &= q^{(\mu+n-\half)/2} [j - \mu]^\half    \ B^-_{jn}.
\end{align*}
The new matrices $A^\pm_{jn}$, $B^\pm_{jn}$ may be further refined by 
using commutation relations involving $\rho'(f)$ and $\rho'(e)$. For 
instance, $\rho'(f) \pi'(a) = q^{-\half} \pi'(a)\,\rho'(f)$ entails
\begin{align*}
&\begin{pmatrix}
[j - n + \half]^\half [j + n + \tfrac{5}{2}]^\half & 0 \\
0 & [j - n - \half]^\half [j + n + \sesq]^\half \end{pmatrix} A^+_{jn}
\\
&\qquad = A^+_{j,n+1} 
\begin{pmatrix} [j - n + \half]^\half [j + n + \sesq]^\half & 0 \\
0 & [j - n - \half]^\half [j + n + \half]^\half \end{pmatrix}.
\end{align*}
This yields four recurrence relations for the entries of $A^+_{jn}$,
one of which has only the trivial solution; we conclude that
$$
A^+_{jn} = \begin{pmatrix} [j + n + \sesq]^\half a^+_{j\up\up}
& 0 \\[\jot]
[j - n + \half]^\half a^+_{j\dn\up} & 
[j + n + \half]^\half a^+_{j\dn\dn} \end{pmatrix},
$$
where the $a^+_{j\updn\updn}$ are scalars depending only on~$j$. In a 
similar fashion, we arrive at
\begin{align*}
A^-_{jn} &= \begin{pmatrix} [j - n + \half]^\half a^-_{j\up\up}
& [j + n + \half]^\half a^-_{j\up\dn} \\[\jot]
0 & [j - n - \half]^\half a^-_{j\dn\dn} \end{pmatrix},
\\[2\jot]
B^+_{jn} &= \begin{pmatrix} [j - n + \sesq]^\half b^+_{j\up\up} & 0
\\[\jot]
[j + n + \half]^\half b^+_{j\dn\up}
& [j - n + \half]^\half b^+_{j\dn\dn} \end{pmatrix},
\\[2\jot]
B^-_{jn} &= \begin{pmatrix} [j + n + \half]^\half b^-_{j\up\up}
& [j - n + \half]^\half b^-_{j\up\dn} \\[\jot]
0 & [j + n - \half]^\half b^-_{j\dn\dn} \end{pmatrix}.
\end{align*}

The analogue of \eqref{eq:ab-link} leads quickly to the relations
\begin{equation}
\begin{aligned}
b^+_{j\up\up} &= q^{j+\half} a^+_{j\up\up},  \\
b^-_{j\up\up} &= - q^{-j-\sesq} a^-_{j\up\up},
\end{aligned}  \qquad
\begin{aligned}
b^+_{j\dn\up} &= - q^{-j-\sesq} a^+_{j\dn\up},  \\
b^-_{j\up\dn} &= q^{j-\half} a^-_{j\up\dn},
\end{aligned}  \qquad
\begin{aligned}
b^+_{j\dn\dn} &= q^{j-\half} a^+_{j\dn\dn},  \\
b^-_{j\dn\dn} &= - q^{-j-\half} a^-_{j\dn\dn}.
\end{aligned}
\label{eq:jonly-ba}
\end{equation}
Next, from the analogue of \eqref{eq:bac-link} we get
$$
(a^-_{j+\half,\up\up})^\star = q^{2j+2} a^+_{j\up\up}, \qquad
(a^-_{j+\half,\up\dn})^\star = - a^+_{j\dn\up}, \qquad
(a^-_{j+\half,\dn\dn})^\star = q^{2j+1} a^+_{j\dn\dn}.
$$

The $a^+_{j\updn\updn}$ parameters may be determined from
$\pi'(b)\pi'(a)\,\kett{j\mu n} = q\,\pi'(a)\pi'(b)\,\kett{j\mu n}$.
The coefficients of $\kett{j \pm 1,\mu + 1,n}$ yield only the relation
\begin{equation}
[2j+1]\, a^+_{j+\half,\dn\dn} a^+_{j\dn\up}
= [2j+3]\, a^+_{j+\half,\dn\up} a^+_{j\up\up}.
\label{eq:ajp-dnup}
\end{equation}
{}From the $\kett{j,\mu + 1,n}$ terms, we obtain
$$
B^-_{j^+n^+} A^+_{jn} + B^+_{j^-n^+} A^-_{jn}
= q^\half (A^-_{j^+n^-} B^+_{jn} + A^+_{j^-n^-} B^-_{jn}).
$$
Comparison of the diagonal entries on both sides gives two more 
relations:
\begin{align*}
{}[2j+1]\, |a^+_{j\dn\up}|^2
&= q^{2j+1} \bigl( [2j+1]\, |a^+_{j-\half,\up\up}|^2
- q [2j+3]\, |a^+_{j\up\up}|^2 \bigr),
\\
{}[2j+1]\, |a^+_{j-\half,\dn\up}|^2
&= q^{2j} \bigl( q [2j+1]\, |a^+_{j\dn\dn}|^2
- [2j-1]\, |a^+_{j-\half,\dn\dn}|^2 \bigr).
\end{align*}
Finally, the expectation of $\pi'(a^*a + q^2 b^*b) = 1$ in the vector
states for $\ket{j\mu n\up}$ and $\ket{j\mu n\dn}$ leads to the
relations
$$
q^{2j} [2j+1]^2 |a^+_{j-\half,\up\up}|^2 = 1,  \qquad
q^{2j} [2j+1]^2 |a^+_{j\dn\dn}|^2 = 1.
$$
Thus all coefficients are now determined, up to a few $j$-dependent 
phases:
\begin{equation}
a^+_{j\up\up} = \zeta_j \frac{q^{-j-\half}}{[2j+2]},  \qquad
a^+_{j\dn\up} = \eta_j \frac{q^\half}{[2j+1]\,[2j+2]},  \qquad
a^+_{j\dn\dn} = \xi_j \frac{q^{-j}}{[2j+1]},
\label{eq:ajp-solved}
\end{equation}
with $|\zeta_j| = |\eta_j| = |\xi_j| = 1$. The relation
\eqref{eq:ajp-dnup} also implies
$\zeta_{j+\half} \eta_j = \eta_{j+\half} \xi_j$. As before, we may 
reset these phases to~$1$ by redefining $\ket{j\mu n\up}$ and
$\ket{j\mu n\dn}$, without breaking the $(\la',\rho')$-equivariance.
Substituting \eqref{eq:ajp-solved} back in previous formulas then 
gives~\eqref{eq:spin-coeff}.
\end{proof}

As already mentioned, formulas \eqref{eq:spin-coeff} for the matrices
$\a_{j\mu n}^\pm$ and $\b_{j\mu n}^\pm$ could have been obtained also
from a direct but tedious computation using equations
\eqref{eq:spinor-basis} and their inverses.

\begin{rem}
\label{rem:prep-goswami}
Were we to consider a representation of $\A$ that need not be
$(\la',\rho')$-equivariant, we could as well have defined our spinor
space, like in \cite{Goswami}, as $\C^2 \ox V$, instead of
$V \ox \C^2$. The Clebsch--Gordan decomposition of $\C^2 \ox V$ would
be that of equation \eqref{eq:CG-decomp}, but the $q$-Clebsch--Gordan
coefficients appearing in \eqref{eq:spinor-basis-dn} and
\eqref{eq:spinor-basis-up} would be different due to the rule for
exchanging the first two columns in $q$-Clebsch--Gordan coefficients
\cite{KlimykS}:
$$
\CGq j l m r s t = \CGq l j m {-s} {-r} {-t} ,
$$
which results in a substitution of $q$ by $q^{-1}$ in
\eqref{eq:spinor-basis-coeffs}. 

However, this is not the correct lifting of the
$(\la,\rho)$-equivariant representation $\pi$ of $\A$ to a
$(\la',\rho')$-equivariant representation of $\A$ on spinor space. We
already noted that $\pi'$ as defined by $\pi \ox \id$ on $V \ox \C^2$
is $(\la',\rho')$-equivariant, directly from $(\la,\rho)$-equivariance
of $\pi$. One checks, simply by working out both sides of equation
\eqref{eq:covar-spin}, that the noncocommutativity of $\U_q(su(2))$
spoils $(\la'',\rho'')$-equivariance of the representation
$\pi'' := \id \ox \pi$ of $\A$ on the tensor product $\C^2 \ox V$,
where we now define $\rho'' := \id \ox \rho$, and
$$
\la''(h) := (\sg_\half \ox \la)(\cop h)
= \sg_\half(\co{h}{1}) \ox \la(\co{h}{2}).
$$
\end{rem}

\section{The equivariant Dirac operator}
\label{sec:Dirac-oper}

Recall the central Casimir element
$C_q = q k^2 + q^{-1} k^{-2} + (q - q^{-1})^2 ef \in \U$. The
symmetric operators $\la'(C_q)$ and $\rho'(C_q)$ on~$\H$, initially
defined with dense domain~$W$, extend to selfadjoint operators
on~$\H$. The finite-dimensional subspaces $W_j^\up$ and $W_j^\dn$ are
their joint eigenspaces:
\begin{align*}
\la'(C_q)\ket{j\mu n\up}  &= (q^{2j+1} + q^{-2j-1}) \,\ket{j\mu n\up},
&
\rho'(C_q)\ket{j\mu n\up} &= (q^{2j+2} + q^{-2j-2}) \,\ket{j\mu n\up},
\\
\la'(C_q)\ket{j\mu n\dn}  &= (q^{2j+1} + q^{-2j-1}) \,\ket{j\mu n\dn},
&
\rho'(C_q)\ket{j\mu n\dn} &= (q^{2j} + q^{-2j}) \,\ket{j\mu n\dn},
\end{align*}
directly from \eqref{eq:spinor-ladder}.

Let $D$ be a selfadjoint operator on $\H$ which commutes strongly with
$\la'(C_q)$ and $\rho'(C_q)$; then the finite-dimensional subspaces
$W_j^\up$ and $W_j^\dn$ reduce $D$. We look for the general form of
such a selfadjoint operator $D$ which is moreover
$(\la',\rho')$-invariant in the sense that it commutes with $\la'(h)$
and $\rho'(h)$, for each $h \in \U_q(su(2))$.

\begin{lem}
\label{lm:eigen-spinors}
The subspaces $W_j^\up$ and $W_j^\dn$ are eigenspaces for $D$.
\end{lem}

\begin{proof}
We may restrict to either the subspace $W_j^\up$ or $W_j^\dn$. Since
$\la'(k)$ and $\rho'(k)$ are required to commute with $D$ and moreover
have distinct eigenvalues on these subspaces, it follows that $D$ has
a diagonal matrix with respect to the basis $\ket{j\mu n\up}$,
respectively $\ket{j\mu n\dn}$. If we provisionally write
$D\ket{j\mu n\up} = d_{j\mu n}^\up \,\ket{j\mu n\up}$, then the
vanishing of
$$
[D,\la'(f)]\,\ket{j\mu n\up} = (d_{j,\mu+1,n}^\up - d_{j\mu n}^\up)\,
[j-\mu]^\half [j+\mu+1]^\half \,\ket{j,\mu+1,n\up},
$$
for $\mu = -j,\dots,j-1$, 
shows that $d_{j\mu n}^\up$ is independent of~$\mu$; and 
$[D,\rho'(f)] = 0$ likewise shows that $d_{j\mu n}^\up$ does not
depend
on~$n$. The same goes for $d_{j\mu n}^\dn$, too. Thus we may write
\begin{equation}
D\ket{j\mu n\up} = d_j^\up \,\ket{j\mu n\up},  \qquad
D\ket{j\mu n\dn} = d_j^\dn \,\ket{j\mu n\dn},
\label{eq:Dirac-eigen}
\end{equation}
where $d_j^\up$ and $d_j^\dn$ are real eigenvalues of~$D$. The
respective multiplicities are $(2j + 1)(2j + 2)$ and $2j(2j + 1)$, in 
view of~\eqref{eq:Wj-dims}.
\end{proof}

One of the conditions for the triple $(\A, \H, D)$ to be a spectral
triple, is boundedness of the commutators $[D,\pi'(x)]$ for
$x \in \A$. This naturally imposes certain restrictions on the
eigenvalues $d_j^\up, d_j^\dn$ of the operator $D$.

For convenience, we recall the representation $\pi'$ of $a$ in the
basis $\kett{j\mu n}$, written explicitly on $\ket{j\mu n\up}$ and
$\ket{j\mu n\dn}$ as in~\eqref{eq:act-on-arrows}:
\begin{align*}
\pi'(a) \ket{j\mu n\up} &= \sum_\pm \a_{j\mu n\up\up}^\pm
\ket{j^\pm\mu^+n^+\up} + \a_{j\mu n\dn\up}^+ \ket{j^+\mu^+n^+\dn},
\\ 
\pi'(a) \ket{j\mu n \dn} &= \sum_\pm \a_{j\mu n\dn\dn}^\pm
\ket{j^\pm\mu^+n^+\dn} + \a_{j\mu n\up\dn}^- \ket{j^-\mu^+n^+\up}.
\end{align*}
Then, a straightforward computation shows that
\begin{align}
[D,\pi'(a)] \,\ket{j\mu n\up} &= \sum_\pm
\a_{j\mu n\up\up}^\pm (d_{j^\pm}^\up - d_j^\up) \ket{j^\pm\mu^+n^+\up}
+ \a_{j\mu n\dn\up}^+ (d_{j^+}^\dn - d_j^\up) \ket{j^+ \mu^+ n^+ \dn},
\nn \\ 
[D,\pi'(a)] \,\ket{j\mu n\dn} &= \sum_\pm
\a_{j\mu n\dn\dn}^\pm (d_{j^\pm}^\dn - d_j^\dn) \ket{j^\pm\mu^+n^+\dn}
+ \a_{j\mu n\up\dn}^- (d_{j^-}^\up - d_j^\dn) \ket{j^- \mu^+ n^+ \up}.
\label{eq:D-commutator}
\end{align} 

Recall that the standard Dirac operator $\Dslash$ on the
sphere $\Sf^3$, with the round metric, has eigenvalues
$(2j + \sesq)$ for $j = 0,\half,1,\sesq$, with respective
multiplicities $(2j + 1)(2j + 2)$; and $-(2j + \half)$ for
$j = \half,1,\sesq$, with respective multiplicities $2j(2j + 1)$: see
\cite{Baer,Homma}, for instance. Notice that its spectrum is symmetric
about~$0$. 

In \cite{BibikovK} a ``$q$-Dirac'' operator $D$ was proposed, which
in our notation corresponds to taking
$d_j^\up = 2[2j+1]/(q + q^{-1})$ and $d_j^\dn = - d_j^\up$; these are
$q$-analogues of the classical eigenvalues of $\Dslash - \half$.
For this particular choice of eigenvalues, it follows directly from
the explicit form \eqref{eq:spin-coeff} of the matrices
$\a_{j\mu n}^\pm$ that then the right hand sides of
\eqref{eq:D-commutator} diverge, and therefore $[D,\pi'(a)]$ is
unbounded. This was already noted in \cite{ConnesL} and it was
suggested that one should instead consider an operator $D$ whose
spectrum matches that of the classical Dirac operator. In fact,
Proposition~\ref{pr:first-order} below shows that this is essentially
the only possibility for a Dirac operator satisfying a (modified)
first-order condition.

Let us then consider any operator $D$ given by \eqref{eq:Dirac-eigen}
--that is, a bi-equivariant one-- with eigenvalues of the following
form:
\begin{equation}
d_j^\up = c_1^\up j + c_2^\up, \qquad
d_j^\dn = c_1^\dn j + c_2^\dn,
\label{eq:linear-evs}
\end{equation}
where $c_1^\up$, $c_2^\up$, $c_1^\dn$, $c_2^\dn$ are independent
of~$j$. For brevity, we shall say that the eigenvalues are ``linear
in~$j$''. On the right hand side of \eqref{eq:D-commutator}, the
``diagonal'' coefficients simplify to 
\begin{equation}
\a_{j\mu n\up\up}^\pm (d_{j^\pm}^\up - d_j^\up)
= \half \a_{j\mu n\up\up}^\pm c_1^\up,  \qquad
\a_{j\mu n\dn\dn}^\pm (d_{j^\pm}^\dn - d_j^\dn)
= \half \a_{j\mu n\dn\dn}^\pm c_1^\dn,
\label{eq:diag-diffs}
\end{equation}
which can be uniformly bounded with respect to~$j$ --see expressions
\eqref{eq:spin-coeff}. For the off-diagonal terms, involving
$\a_{j\mu n \dn\up}^+$ and $\a_{j\mu n\up\dn}^-$, the differences
between the ``up'' and ``down'' eigenvalues are linear in~$j$. Since
$0 < q < 1$, it is clear that $[N] \sim (q^{-1})^{N-1}$ for large~$N$,
and thus $\a_{j\mu n\dn\up}^+ \sim q^{3j+n+\sesq} \leq q^{2j+1}$ for
large~$j$. Similar easy estimates yield
\begin{equation}
\begin{aligned}
\a_{j\mu n\dn\up}^+ &= O(q^{2j+1}), 
\\
\a_{j\mu n\up\dn}^- &= O(q^{2j}),
\end{aligned}
\qquad
\begin{aligned}
\b_{j\mu n\dn\up}^+ &= O(q^{2j+\half}), 
\\
\b_{j\mu n\up\dn}^- &= O(q^{2j+\half}), \quad\text{as } j \to \infty.
\end{aligned}
\label{eq:est-offdiag}
\end{equation} 
We therefore arrive at
\begin{equation}
|\a_{j\mu n\dn\up}^+ (d_{j^+}^\dn - d_j^\up - 1)| \leq C j q^{2j},
\qquad
|\a_{j\mu n\up\dn}^- (d_{j^-}^\up - d_j^\dn - 1)| \leq C' j q^{2j},
\label{eq:Dcomm-off-diag}
\end{equation} 
for some $C > 0$, $C' > 0$, independent of~$j$; and similar estimates 
hold for the off-diagonal coefficients of $\pi'(b)$.

\begin{prop}
\label{pr:D-comm}
Let $D$ be any selfadjoint operator with eigenspaces $W_j^\up$ and
$W_j^\dn$, and eigenvalues \eqref{eq:Dirac-eigen}. If the eigenvalues
$d_j^\up$ and $d_j^\dn$ are linear in~$j$ as in \eqref{eq:linear-evs},
then $[D,\pi'(x)]$ is a bounded operator for all $x \in \A$.
\end{prop}

\begin{proof}
Since $a$ and $b$ generate~$\A$ as a $*$-algebra, it is enough to
consider the cases $x = a$ and $x = b$. For $x = a$ and any
$\xi \in \H$, the relations \eqref{eq:D-commutator} and
\eqref{eq:diag-diffs}, together with the Schwarz inequality, give the
estimate
$$
\|[D,\pi'(a)]\,\xi\|^2
\leq \tfrac{1}{4} \max\{(c_1^\up)^2, (c_1^\dn)^2\}
\,\|\pi'(a)\xi\|^2 + \|\xi\|^2 \|\eta\|^2 ,
$$
where $\eta$ is a vector whose components are estimated by
\eqref{eq:Dcomm-off-diag}, which establishes finiteness of~$\|\eta\|$
since $0 < q < 1$. Therefore, $[D,\pi'(a)]$ is norm bounded. In the
same way, we find that $[D,\pi'(b)]$ is bounded.
\end{proof}

Now, if $D$ is a selfadjoint operator as in
Proposition~\ref{pr:D-comm}, and if the eigenvalues of $D$ 
satisfy \eqref{eq:linear-evs} and, moreover,
\begin{equation}
c_1^\dn = - c_1^\up,  \qquad  c_2^\dn = - c_2^\up + c_1^\up,
\label{eq:nice-evs}
\end{equation}
then the spectrum of~$D$ coincides with that of the classical
Dirac operator $\Dslash$ on the round sphere $\Sf^3$, up to rescaling
and addition of a constant. Thus, we can regard our spectral triple
as an isospectral deformation of $(C^\infty(\Sf^3),\H,\Dslash)$, and
in particular, its spectral dimension is~$3$. We summarize our
conclusions in the following theorem.

\begin{thm}
The triple $(\A(SU_q(2)),\H,D)$, where the eigenvalues of~$D$ satisfy
\eqref{eq:linear-evs} and \eqref{eq:nice-evs}, is a $3^+$-summable
spectral triple.
\qed
\end{thm}

At this point, it is appropriate to comment on the relation of our
construction with that of~\cite{Goswami}. There, a spinor
representation is constructed by tensoring the left regular
representation of $\A(SU_q(2))$ by $\C^2$ on the left. This spinor
space is then decomposed into two subspaces, similar to our ``up'' and
``down'' subspaces, on which $D$ acts diagonally with eigenvalues
linear in the total spin number~$j$. The corresponding decomposition
of the representation $\pi'$ of $\A(SU_q(2))$ on spinor space is
obtained by using the appropriate Clebsch--Gordan coefficients.
However, contrary to what we have established above, in \cite{Goswami}
it is found that a certain commutator $[D,\pi'(x)]$ is an unbounded
operator. In particular, the off-diagonal terms in the representation of \cite{Goswami}
do not have the compact nature we encountered in
\eqref{eq:est-offdiag}. They can be bounded from below by a positive
constant, which leads, when multiplied by a term linear in~$j$, to an
unbounded operator.

The origin of this notable contrast is the following. Since in
\cite{Goswami} no condition of $\U_q(su(2))$-equivariance is imposed
\textit{a priori} on the representation of $\A(SU_q(2))$, the spinor
space $W$ could be identified either with $V \ox \C^2$ or
$\C^2 \ox V$, according to convenience. However, as we noted in Remark
\ref{rem:prep-goswami}, the choice of $\C^2 \ox V$ is not allowed by
the condition of $(\la',\rho')$-equivariance, because $\U_q(su(2))$ is
not cocommutative. Indeed, repeating the construction of a spinor
representation and Dirac operator on the spinor space
\mbox{$\C^2 \ox V$} instead of $V \ox \C^2$ --hence ignoring
equivariance-- results eventually in unbounded commutators.

\section{The real structure}
\label{sec:J-oper}

The next issue we address is the real structure $J$ on the spectral
triple $(\A(SU_q(2)),\H,D)$. We shall see that by requiring
equivariance of $J$ it is not possible to satisfy all usual properties
of a real spectral triple like in \cite{ConnesReal} or \cite{Polaris}.
Among other things, these conditions entail for $J$ that it intertwine
a left action and a commuting right action of the algebra on the
Hilbert space, which then gets a bimodule structure (the commutant
property); and that the bounded commutators $[D,a]$, for any element
$a$ in the algebra, commute with the opposite action by any $b$ in the
algebra (the first order condition on~$D$). However, we shall be able
to satisfy these two conditions only up to certain compact operators.

\subsection{The Tomita operator of the regular representation}
\label{ssc:Tomita-oper}

On the GNS representation space $\H_\psi$, there is a natural
involution $T_\psi \: \eta(x) \mapsto \eta(x^*)$, with domain
$\eta(C(SU_q(2)))$, which may be regarded as an unbounded (antilinear)
operator on~$\H_\psi$.
The Tomita--Takesaki theory \cite{Takesaki} shows that this operator
is closable (we denote its closure also by~$T_\psi$) and that the
polar decomposition $T_\psi =: J_\psi \Delta_\psi^{1/2}$ defines both
the positive ``modular operator'' $\Delta_\psi$ and the antiunitary
``modular conjugation'' $J_\psi$. It has already been noted by
Chakraborty and Pal~\cite{ChakrabortyPDual} that this $J_\psi$ has a
simple expression in terms of the matrix elements of our chosen 
orthonormal basis for~$\H_\psi$. Indeed, it follows immediately 
from~\eqref{eq:matelt-star} and \eqref{eq:matelt-onb} that 
$$
T_\psi \,\ket{lmn} = (-1)^{2l+m+n} q^{m+n} \,\ket{l,-m,-n}.
$$
One checks, using \eqref{eq:reg-repn}, that 
$$
T_\psi \pi(a)\,\ket{000} = \pi(a^*) \,\ket{000},  \qquad
T_\psi \pi(b)\,\ket{000} = \pi(b^*) \,\ket{000}.
$$
Since $\pi$ is the GNS representation for the state~$\psi$, this is 
enough to conclude that
\begin{equation}
T_\psi \eta(x) = \eta(x^*)  \sepword{for all}  x \in \A.
\label{eq:Tomita-prop}
\end{equation}

The adjoint antilinear operator, satisfying
$\braCket{\eta}{T^*_\psi}{\xi} = \braCket{\xi}{T_\psi}{\eta}$, is 
given by
$T^*_\psi \,\ket{lmn} = (-1)^{2l+m+n} q^{-m-n} \,\ket{l,-m,-n}$, and
since $\Delta_\psi = T^*_\psi T_\psi$, we see that every $\ket{lmn}$
lies in $\Dom\Delta_\psi$ with
$\Delta_\psi \,\ket{lmn} = q^{2m+2n} \,\ket{lmn}$. Consequently,
\begin{equation}
J_\psi \,\ket{lmn} = (-1)^{2l+m+n} \,\ket{l,-m,-n}.
\label{eq:Jpsi-formula}
\end{equation}
It is clear that $J_\psi^2 = 1$ on~$\H_\psi$.

\begin{defn}
\label{df:piop-repn}
Let $\piop(x) := J_\psi\, \pi(x^*) \,J_\psi^{-1}$, so that $\piop$ is
a $*$-antirepresentation of $\A$ on~$\H_\psi$. Equivalently, $\piop$
is a $*$-representation of the opposite algebra $\A(SU_{1/q}(2))$.
By Tomita's theorem~\cite{Takesaki}, $\pi$ and $\piop$ are commuting
representations.
\end{defn}

As an example, we compute
\begin{align*}
\piop(a)\,\ket{lmn}
&= (-1)^{2l+m+n} J_\psi \pi(a^*) \,\ket{l,-m,-n}
\\
&= (-1)^{2l+m+n} J_\psi \bigl( \Onda A^+_{l,-m,-n} \ket{l^+,-m^+,-n^+}
 + \Onda A^-_{l,-m,-n} \ket{l^-,-m^+,-n^+} \bigr)
\\
&= \Onda A^+_{l,-m,-n} \ket{l^+ m^+ n^+}
 + \Onda A^-_{l,-m,-n} \ket{l^- m^+ n^+}
\\
&= A^-_{l^+,-m^+,-n^+} \ket{l^+ m^+ n^+}
 + A^+_{l^-,-m^+,-n^+} \ket{l^- m^+ n^+},
\end{align*}
where, explicitly,
\begin{align*}
A^-_{l^+,-m^+,-n^+} &= q^{(2l-m-n+1)/2}
\biggl( \frac{[l+m+1][l+n+1]}{[2l+1][2l+2]} \biggr)^\half \!,
\\
A^+_{l^-,-m^+,-n^+} &= q^{-(2l+m+n+1)/2} 
\biggl( \frac{[l-m][l-n]}{[2l][2l+1]} \biggr)^\half \!.
\end{align*}
A glance back at \eqref{eq:reg-repn-coeffs} shows that these
coefficients are identical with those of $\pi(a)\,\ket{lmn}$,
\textit{after substituting $q \mapsto q^{-1}$}. A similar phenomenon
occurs with the coefficients of $\piop(b)$. We find, indeed, that
\begin{align*}
\piop(a) \,\ket{lmn}
&= A^{\circ+}_{lmn} \ket{l^+ m^+ n^+}
 + A^{\circ-}_{lmn} \ket{l^- m^+ n^+},
\\
\piop(b) \,\ket{lmn}
&= B^{\circ+}_{lmn} \ket{l^+ m^+ n^-}
 + B^{\circ-}_{lmn} \ket{l^- m^+ n^-},
\end{align*}
where
\begin{equation}
A^{\circ\pm}_{lmn}(q) = A^\pm_{lmn}(q^{-1}),  \qquad
B^{\circ\pm}_{lmn}(q) = q^{-1} B^\pm_{lmn}(q^{-1}).
\label{eq:piop-elts}
\end{equation}

We can now verify directly that the representations $\pi$ and $\piop$ 
commute, without need to appeal to the theorem of Tomita. For 
instance,
\begin{gather*}
\braCket{l+1,m+1,n+1}{[\pi(a),\piop(a)]}{lmn}
= A^{\circ+}_{l^+m^+n^+} A^+_{lmn} - A^+_{l^+m^+n^+} A^{\circ+}_{lmn}
\\
= Q \biggl(
\frac{[l+m+1][l+m+2][l+n+1][l+n+2]}{[2l+1][2l+2]^2[2l+3]}
\biggr)^\half \!,
\end{gather*}
where 
$$
Q =  q^{\half(2l^+ - m^+ - n^+ + 1)} q^{\half(-2l + m + n - 1)}
  - q^{\half(-2l^+ + m^+ + n^+ - 1)} q^{\half( 2l - m - n + 1)} = 0.
$$
Likewise, $\braCket{l-1,m+1,n+1}{[\pi(a),\piop(a)]}{lmn} = 0$, and
one checks that the matrix element
$\braCket{l,m+1,n+1}{[\pi(a),\piop(a)]}{lmn}$ vanishes, too.

\vspace{6pt}

The $(\la,\rho)$-equivariance of~$\pi$ is reflected in an analogous 
equivariance condition for~$\piop$. We now identify this condition 
explicitly.

\begin{lem}
\label{lm:piop-eqvar}
The symmetry of the antirepresentation $\piop$ of $\A$ on~$\H_\psi$ 
is given by the equivariance conditions:
\begin{align}
\la(h)\,\piop(x) \xi
&= \piop(\co{\tilde h}{2} \. x) \, \la(\co{h}{1}) \xi,
\nn \\[\jot]
\rho(h)\,\piop(x)\xi
&= \piop(\co{\tilde h}{2} \lt x)\,\rho(\co{h}{1}) \xi,
\label{eq:covar-piop}
\end{align}
for all $h \in \U$, $x \in \A$ and $\xi \in V$, and
$h \mapsto \tilde h$ is the automorphism of~$\U$ determined on 
generators by $\tilde k := k$, $\tilde f := q^{-1} f$, and
$\tilde e := q e$.
\end{lem}

\begin{proof}
We work only on the dense subspace $V$. From \eqref{eq:uqsu2-repns}
and~\eqref{eq:Jpsi-formula}, we get at once
\begin{equation}
J_\psi \la(k)^* J_\psi^{-1} = \la(k^{-1}),  \quad
J_\psi \la(f)^* J_\psi^{-1} = - \la(f),  \quad
J_\psi \la(e)^* J_\psi^{-1} = - \la(e),
\label{eq:la-conj}
\end{equation}
and identical relations with $\rho$ instead of~$\la$. Write $\a$ for 
the antiautomorphism of~$\U$ determined by $\a(k) := k^{-1}$,
$\a(f) := -f$, and $\a(e) := -e$; so that 
$J_\psi \la(h)^* J_\psi^{-1} = \la(\a(h))$ for $h \in \U$, and 
similarly with $\rho$ instead of~$\la$. 

Next, the first relation of~\eqref{eq:covar-repn} is equivalent to
\begin{equation}
\pi(x)\,\la(Sh) = \la(S\co{h}{1}) \, \pi(\co{h}{2} \. x).
\label{eq:eqvar-switch}
\end{equation}
Indeed, the left hand side can be expanded as
$$
\pi(x)\,\la(\eps(\co{h}{1})\,S\co{h}{2}) 
= \la(S\co{h}{1}\,\co{h}{2})\,\pi(x)\,\la(S\co{h}{3}) 
= \la(S\co{h}{1})\,\pi(\co{h}{2}\.x)\,\la(\co{h}{3})\,\la(S\co{h}{4})
$$
on applying \eqref{eq:covar-repn}; and the rightmost expression
equals the right hand side of~\eqref{eq:eqvar-switch}. Taking 
hermitian adjoints and conjugating by~$J_\psi$, we get 
$$
\la(\a(Sh))\,\piop(x) = \piop(\co{h}{2} \. x) \, \la(\a(S\co{h}{1})).
$$
It remains only to note that $S\a = \a S$ is an automorphism of~$\U$, 
whose inverse is the map $h \mapsto \tilde h$ above; and to repeat
the 
argument with $\rho$ instead of~$\la$, changing only the left action 
of~$\U$ in concordance with~\eqref{eq:covar-repn}.
\end{proof}

\vspace{6pt}

An independent check of \eqref{eq:covar-piop} is afforded by the
following argument. We may ask which antirepresentations $\piop$
of~$\H_\psi$ satisfy these equivariance conditions. It suffices to run
the proof of Proposition~\ref{pr:suq2-repn}, \emph{mutatis mutandis},
to determine the possible form of such a $\piop$ on the basis vectors
$\ket{lmn}$. For instance, \eqref{eq:ab-link} is replaced by
$$
\rho(e) \piop(a)
= \piop(\tilde e \lt a) \rho(k^{-1}) + \piop(\tilde k \lt a) \rho(e)
= q\,\piop(b) \rho(k^{-1}) + q^\half \piop(a) \rho(e).
$$
One finds that all formulas in that proof are reproduced, except for 
changes in the powers of~$q$ that appear; and, apart from the 
aforementioned phase ambiguities, one recovers precisely the form of 
$\piop$ given by \eqref{eq:piop-elts}.

\vspace{6pt}

Before proceeding, we indicate also the symmetry of the Tomita
operator $T_\psi$, analogous to \eqref{eq:la-conj} above. Combining
\eqref{eq:Tomita-prop} with \eqref{eq:covar-repn}, and recalling that
$\eta(x) = \pi(x)\,\ket{000}$, we find that for generators $h$ of
$\U$,
$$
T_\psi \la(h) \pi(x) \,\ket{000} = \pi(x^* \rt \vth(h)^*)\,\ket{000}.
$$
On the other hand,
$$
\la (\vth^{-1} S( \vth(h^*))) T_\psi \pi(x)\,\ket{000}
= \pi(x^* \rt \vth(h)^*) \,\ket{000}.
$$
One checks easily on generators that
$\vth^{-1} S(\vth(h)^*) = S(h)^*$. Since the vector $\ket{000}$ is 
separating for the GNS representation, we conclude that 
$$
T_\psi \,\la(h)\, T_\psi^{-1} = \la(Sh)^* .
$$
Similarly, we find that
$$
T_\psi \,\rho(h)\, T_\psi^{-1} = \rho(Sh)^* .
$$
In other words, the antilinear involutory automorphism
$h \mapsto (Sh)^*$ of the Hopf $*$-algebra~$\U$ is implemented by the 
Tomita operator for the Haar state of the dual Hopf $*$-algebra $\A$.
This is a known feature of quantum-group duality in the $C^*$-algebra 
framework; for this and several other implementations by spatial 
operators, see~\cite{MasudaNW}.

\subsection{The real structure on spinors}
\label{ssc:spinor-J}

We are now ready to come back to spinors. Notice that $J_\psi$ does
not appear explicitly in the equivariance conditions
\eqref{eq:covar-piop} for the right regular representation~$\piop$ of
$\A$ on~$\H_\psi$. Thus, we are now able to construct the ``right
multiplication'' representation of~$\A$ on spinors from its symmetry
alone, and to deduce the conjugation operator $J$ on spinors after the
fact.

\begin{prop}
\label{pr:piiop-eqvar}
Let $\piiop$ be an antirepresentation of $\A$ on 
$\H = \H_\psi \oplus \H_\psi$ satisfying the following equivariance
conditions:
\begin{align}
\la'(h)\,\piiop(x) \xi
&= \piiop(\co{\tilde h}{2} \. x) \, \la'(\co{h}{1}) \xi,
\nn \\[\jot]
\rho'(h)\,\piiop(x)\xi
&= \piiop(\co{\tilde h}{2} \lt x)\,\rho'(\co{h}{1}) \xi.
\label{eq:covar-piiop}
\end{align}
Then, up to some phase factors depending only on the index~$j$ in the 
decomposition \eqref{eq:spinor-decomp}, $\piiop$ is given on the 
spinor basis by
\begin{align}
\piiop(a) \,\kett{j\mu n}
&= \a^{\circ+}_{j\mu n} \kett{j^+ \mu^+ n^+}
 + \a^{\circ-}_{j\mu n} \kett{j^- \mu^+ n^+},
\nn \\[\jot]
\piiop(b) \,\kett{j\mu n}
&= \b^{\circ+}_{j\mu n} \kett{j^+ \mu^+ n^-}
 + \b^{\circ-}_{j\mu n} \kett{j^- \mu^+ n^-},
\nn \\[\jot]
\piiop(a^*) \,\kett{j\mu n}
&= \tilde\a^{\circ+}_{j\mu n} \kett{j^+ \mu^- n^-}
 + \tilde\a^{\circ-}_{j\mu n} \kett{j^- \mu^- n^-},
\label{eq:right-spin-repn}
\\[\jot]
\piiop(b^*) \,\kett{j\mu n}
&= \tilde\b^{\circ+}_{j\mu n} \kett{j^+ \mu^- n^+}
 + \tilde\b^{\circ-}_{j\mu n} \kett{j^- \mu^- n^+},
\nn
\end{align}
where $\a^{\circ\pm}_{j\mu n}$ and $\b^{\circ\pm}_{j\mu n}$ are the
triangular $2 \x 2$ matrices, given by
$\a^{\circ\pm}_{j\mu n}(q) = \a^\pm_{j\mu n}(q^{-1})$ and
$\b^{\circ\pm}_{j\mu n}(q) = q^{-1} \b^\pm_{j\mu n}(q^{-1})$, with
$\a^\pm_{j\mu n}$ and $\b^\pm_{j\mu n}$ given
by~\eqref{eq:spin-coeff}.
\end{prop}

\begin{proof}
We retrace the steps of the proof of Proposition~\ref{pr:spin-repn},
\emph{mutatis mutandis}. Since $\tilde k \. a = k \. a = q^\half a$,
the relations involving $\la'(k)$ and $\rho'(k)$ are unchanged. We 
quickly conclude that $\piiop$ must have the form 
\eqref{eq:right-spin-repn}, and it remains to determine the 
coefficient matrices.

The commutation relations of $\la'(f)$ with $\piiop(a)$ and
$\piiop(b)$ give:
\begin{align*}
\a^{\circ+}_{j\mu n}
&= q^{-\half(\mu+n-\half)} [j + \mu + 1]^\half\ A^{\circ+}_{jn}, &
\a^{\circ-}_{j\mu n}
&= q^{-\half(\mu+n-\half)} [j - \mu]^\half    \ A^{\circ-}_{jn}, \\
\b^{\circ+}_{j\mu n}
&= q^{-\half(\mu+n-\half)} [j + \mu + 1]^\half\ B^{\circ+}_{jn}, &
\b^{\circ-}_{j\mu n}
&= q^{-\half(\mu+n-\half)} [j - \mu]^\half    \ B^{\circ-}_{jn}.
\end{align*}

The matrices $A^{\circ\pm}_{jn}$, $B^{\circ\pm}_{jn}$ may be
determined, as before, by the commutation relations involving
$\rho'(f)$ and $\rho'(e)$. One finds that the $n$-dependent factors
such as $[j + n + \sesq]^\half$ and so on, are the same as the 
respective entries of $A^\pm_{jn}$, $B^\pm_{jn}$; let 
$a^{\circ+}_{j\up\up}$, etc., be the remaining factors which depend 
on~$j$ only. Then \eqref{eq:jonly-ba} is replaced by
\begin{align*}
b^{\circ+}_{j\up\up} &= q^{-j-\sesq} a^{\circ+}_{j\up\up}, &
b^{\circ+}_{j\dn\up} &= - q^{j+\half} a^{\circ+}_{j\dn\up}, &
b^{\circ+}_{j\dn\dn} &= q^{-j-\half} a^{\circ+}_{j\dn\dn},
\\
b^{\circ-}_{j\up\up} &= - q^{j+\half} a^{\circ-}_{j\up\up}, &
b^{\circ-}_{j\up\dn} &= q^{-j-\half} a^{\circ-}_{j\up\dn}, &
b^{\circ-}_{j\dn\dn} &= - q^{j-\half} a^{\circ-}_{j\dn\dn}.
\end{align*}
Next, we find
$$
(a^{\circ-}_{j+\half,\up\up})^\star = q^{-2j-2} a^{\circ+}_{j\up\up},
\qquad
(a^{\circ-}_{j+\half,\up\dn})^\star = - a^{\circ+}_{j\dn\up},
\qquad
(a^{\circ-}_{j+\half,\dn\dn})^\star = q^{-2j-1} a^{\circ+}_{j\dn\dn}.
$$

Since $\piiop$ is an antirepresentation, $ab = q^{-1} ba$
implies $\piiop(b)\piiop(a) = q^{-1}\,\piiop(a)\piiop(b)$. The matrix 
elements of both sides lead to three relations:
\begin{equation}
[2j+1]\, a^{\circ+}_{j+\half,\dn\dn} a^{\circ+}_{j\dn\up}
= [2j+3]\, a^{\circ+}_{j+\half,\dn\up} a^{\circ+}_{j\up\up},
\label{eq:ajp-dnup-opp}
\end{equation}
which is formally identical to \eqref{eq:ajp-dnup}, and
\begin{align*}
[2j+1]\, |a^{\circ+}_{j\dn\up}|^2
&= q^{-2j-1} \bigl( [2j+1]\, |a^{\circ+}_{j-\half,\up\up}|^2
- q^{-1} [2j+3]\, |a^{\circ+}_{j\up\up}|^2 \bigr),
\\
[2j+1]\, |a^{\circ+}_{j-\half,\dn\up}|^2
&= q^{-2j} \bigl( q^{-1} [2j+1]\, |a^{\circ+}_{j\dn\dn}|^2
- [2j-1]\, |a^{\circ+}_{j-\half,\dn\dn}|^2 \bigr).
\end{align*}

Finally, the relation $aa^* + bb^* = 1$ yields
$\piiop(a^*)\piiop(a) + \piiop(b^*)\piiop(b) = 1$; its diagonal
matrix elements gives the last two relations:
$$
q^{-2j} [2j+1]^2 |a^{\circ+}_{j-\half,\up\up}|^2 = 1,  \qquad
q^{-2j} [2j+1]^2 |a^{\circ+}_{j\dn\dn}|^2 = 1.
$$
All coefficients are now determined except for their phases:
\begin{equation}
a^{\circ+}_{j\up\up} = \zeta_j^\circ \frac{q^{j+\half}}{[2j+2]},
\quad
a^{\circ+}_{j\dn\up} = \eta_j^\circ \frac{q^{-\half}}{[2j+1]\,[2j+2]},
\quad
a^{\circ+}_{j\dn\dn} = \xi_j^\circ \frac{q^j}{[2j+1]},
\label{eq:ajp-solved-opp}
\end{equation}
and \eqref{eq:ajp-dnup-opp} also entails the phase relations
$\zeta_j^\circ \eta_{j+\half}^\circ =
\eta_j^\circ \xi_{j+\half}^\circ$. Once more, we choose all phases to
be~$+1$ by convention. Substituting \eqref{eq:ajp-solved-opp} back in
previous formulas, we find
\begin{equation}
\a^{\circ\pm}_{j\mu n}(q) = \a^\pm_{j\mu n}(q^{-1}),  \qquad
\b^{\circ\pm}_{j\mu n}(q) = q^{-1} \b^\pm_{j\mu n}(q^{-1}).
\label{eq:piiop-elts}
\end{equation}
in perfect analogy with \eqref{eq:piop-elts}.
\end{proof}

\begin{defn}
The conjugation operator $J$ is the antilinear operator on~$\H$ which
is defined explicitly on the orthonormal spinor basis by
\begin{align}
J\, \ket{j\mu n\up} &:= i^{2(2j+\mu+n)} \,\ket{j,-\mu,-n,\up},
\nn \\
J\, \ket{j\mu n\dn} &:= i^{2(2j-\mu-n)} \,\ket{j,-\mu,-n,\dn}.
\label{eq:J-formula}
\end{align}
It is immediate from this presentation that $J$ is antiunitary and
that $J^2 = -1$, since each $4j \pm 2(\mu + n)$ is an odd integer. 
\end{defn}

\begin{prop}
\label{pr:DJ-comm}
The invariant operator $D$ of Section~\ref{sec:Dirac-oper}
commutes with the conjugation operator $J$:
\begin{equation}
J D J^{-1} = D.
\label{eq:DJ-comm}
\end{equation}
\end{prop}

\begin{proof}
This is clear from the diagonal form of both $D$ and $J$ on their 
common eigenspaces $W_j^\up$ and $W_j^\dn$, given by the respective
equations \eqref{eq:Dirac-eigen} and \eqref{eq:J-formula}.
\end{proof}

\begin{rem}
Proposition~\ref{pr:DJ-comm} is a minimal requirement for
$(\A(SU_q(2)),\H,D,J)$ to constitute a real spectral triple. However,
here is where we part company with the axiom scheme for real spectral
triples proposed in~\cite{ConnesReal}. Indeed, the conjugation
operator $J$ that we have defined by~\eqref{eq:J-formula} is
\textit{not} the modular conjugation for the spinor representation of
$\A$. That modular operator is $J_\psi \oplus J_\psi$, which does not
have a diagonal form in our chosen spinor basis (unless $q = 1$). It
is clear that conjugation of $\pi'(\A(SU_q(2))$ by the modular
operator would yield a representation of the opposite algebra
$\A(SU_{1/q}(2))$, and the commutation relation analogous to
\eqref{eq:DJ-comm} would then force $D$ to be equivariant under the
corresponding symmetry of $U_{1/q}(su(2))$, denoted by
$(\la'',\rho'')$ in our earlier Remark~\ref{rem:prep-goswami}. It is
not hard to check that this extra equivariance condition would force
$D$ to be merely a scalar operator, thereby negating the possibility
of an equivariant $3^+$-summable real spectral triple based on
$\A(SU_q(2))$ with the modular conjugation operator. This result is
consonant with the ``no-go theorem'' of Schm\"udgen \cite{Schmuedgen} 
for nontrivial commutator representations of Woronowicz differential
calculi on $SU_q(2)$.

The remedy that we propose here is to modify $J$, in keeping with the 
symmetry of the spinor representation, to a non-Tomita conjugation
operator. We shall see, however, that the expected properties of
real spectral triples do hold ``up to compact perturbations''.
\end{rem}

It should be noted that $J$ satisfies the analogue
of~\eqref{eq:la-conj} for the representations $\la'$ and~$\rho'$:
\begin{equation}
\begin{aligned}
J \la'(k) J^{-1}  &= \la'(k^{-1}),
\\
J \rho'(k) J^{-1} &= \rho'(k^{-1}),
\end{aligned}
\qquad
\begin{aligned}
J \la'(e) J^{-1} &= -\la'(f),
\\
J \rho'(e) J^{-1} &= -\rho'(f),
\end{aligned}
\label{eq:J-symmetry}
\end{equation}
which follows directly from the definition \eqref{eq:J-formula} and
the relations \eqref{eq:spinor-ladder}.

\begin{prop}
The antiunitary operator $J$ intertwines the left and right spinor
representations:
\begin{equation}
J\, \pi'(x^*) \,J^{-1} = \piiop(x),  \sepword{for all} x \in \A.
\label{eq:J-switch}
\end{equation}
\end{prop}

\begin{proof}
It follows directly from the proof of Lemma~\ref{lm:piop-eqvar},
using the relations \eqref{eq:J-symmetry} instead of 
\eqref{eq:la-conj}, that the antirepresentation 
$x \mapsto J\,\pi'(x^*)\,J^{-1}$ complies with the equivariance
conditions \eqref{eq:covar-piiop}. By 
Proposition~\ref{pr:piiop-eqvar}, it coincides with $\piiop$ up to an 
equivalence obtained by resetting the phase factors in 
\eqref{eq:ajp-solved-opp}. It remains only to check that
$\zeta_j^\circ = \eta_j^\circ = \xi_j^\circ = 1$ for the 
aforementioned antirepresentation. This check is easily effected by
calculating $J\,\pi'(a^*)\,J^{-1}$ directly on the basis vectors
$\ket{j\mu n\up}$. We compute
\begin{align*}
& J \pi'(a^*) J^{-1} \,\ket{j\mu n\up}
= i^{-2(2j-\mu-n)} J \pi'(a^*) \,\ket{j,-\mu,-n,\up}
\\
&= i^{-2(2j-\mu-n)} J \bigl(
\tilde\a^+_{j,-\mu,-n,\up\up} \,\ket{j^+,-\mu^+,-n^+\up}
 + \tilde\a^+_{j,-\mu,-n,\dn\up} \,\ket{j^+,-\mu^+,-n^+\dn}
\\
&\qquad + \tilde\a^-_{j,-\mu,-n,\up\up} \,\ket{j^-,-\mu^+,-n^+\up}
\bigr)
\\
&= \tilde\a^+_{j,-\mu,-n,\up\up} \,\ket{j^+ \mu^+ n^+\up}
 - \tilde\a^+_{j,-\mu,-n,\dn\up} \,\ket{j^+ \mu^+ n^+\dn}
 + \tilde\a^-_{j,-\mu,-n,\up\up} \,\ket{j^- \mu^+ n^+\up}
\\
&= \a^-_{j^+,-\mu^+,-n^+,\up\up} \,\ket{j^+ \mu^+ n^+\up}
 - \a^-_{j^+,-\mu^+,-n^+,\dn\up} \,\ket{j^+ \mu^+ n^+\dn}
 + \a^+_{j^-,-\mu^+,-n^+,\up\up} \,\ket{j^- \mu^+ n^+\up}
\\
&= q^{-\half(\mu+n-\half)} \biggl( 
   q^{j+\half} \frac{[j+\mu+1]^\half [j+n+\sesq]^\half}{[2j+2]}
   \,\ket{j^+ \mu^+ n^+\up}
\\ &\qquad
 + q^{-\half} \frac{[j+\mu+1]^\half [j-n+\half]^\half}{[2j+1][2j+2]}
   \,\ket{j^+ \mu^+ n^+\dn}
 + q^{-j-1} \frac{[j-\mu]^\half [j-n+\half]^\half}{[2j+1]}
   \,\ket{j^- \mu^+ n^+\up} \biggr)
\\
&= \a^{\circ+}_{j\mu n\up\up} \,\ket{j^+ \mu^+ n^+\up}
 + \a^{\circ+}_{j\mu n\dn\up} \,\ket{j^+ \mu^+ n^+\dn}
 + \a^{\circ-}_{j\mu n\up\up} \,\ket{j^- \mu^+ n^+\up}
\\
&= \piiop(a) \,\ket{j\mu n\up},
\end{align*}
where the $\a^{\circ\pm}_{j\mu n}$ coefficients are taken according to
\eqref{eq:piiop-elts}.

In the same way, one finds that
$J \pi'(b^*) J^{-1} \,\ket{j\mu n\up} = \piiop(b) \,\ket{j\mu n\up}$,
again using \eqref{eq:piiop-elts} for $\b^{\circ\pm}_{j\mu n}$; and
similar calculations show that both sides of \eqref{eq:J-switch}
coincide on the basis vector $\ket{j\mu n\dn}$. (These four
calculations, taken together, afford a direct proof of
\eqref{eq:J-switch} without need to consider the symmetries of~$J$.)
\end{proof}

\section{Algebraic properties of the spectral triple}
\label{sec:spec-tri}

In this section, we discuss the properties of the real spectral
triple $(\A(SU_q(2)),\H, D,J)$, in particular its commutant property
and its first-order condition. We will see that these are only
satisfied up to certain compact operators, quite similarly
to~\cite{DabrowskiLPS}.

We can simplify our discussion somewhat by replacing the spinor
representation $\pi'$ of $\A = \A(SU_q(2))$ of
Proposition~\ref{pr:spin-repn} by a so-called approximate
representation $\piappr \: \A \to \B(\H)$, such that
$\pi'(x) - \piappr(x)$ is a compact operator for each
$x \in \A$. In other words, although $\piappr$ need not
preserve the algebra relations of~$\A$, the mappings $\pi'$ and
$\piappr$ have the same image in the Calkin algebra $\B(\H)/\K(\H)$,
that is, they define the same $*$-homomorphism of~$\A$ into the Calkin
algebra.

We denote by $L_q$ the positive trace-class operator given by
$$
L_q \kett{j \mu n} := q^j \,\kett{j \mu n}
\sepword{for} j \in \half\N,
$$
and let $\K_q$ be the two-sided ideal of $\B(\H)$ generated by $L_q$; it
consists of trace-class operators. The ideal $\K_q$ is indeed contained in the ideal of infinitesimals
of order $\a$, that is, compact operators whose $n$-th singular value
$\mu_n$ satisfies $\mu_n = O(n^{-\a})$, for all $\a > 0$. Thus the
following analysis holds modulo infinitesimals of arbitrary high order.

\begin{prop}
\label{pr:appr-rep}
The following equations define  a mapping
$\piappr \: \A \to \B(\H)$ on generators, which is a
$*$-representation modulo $\K_q$, and is approximate to the spin
representation $\pi'$ of Proposition~\ref{pr:spin-repn} in the sense
that $\pi'(x) - \piappr(x) \in \K_q$ for each $x \in \A$:
\begin{align}
\piappr(a) \,\kett{j\mu n}
&= \aappr^+_{j\mu n} \kett{j^+ \mu^+ n^+}
 + \aappr^-_{j\mu n} \kett{j^- \mu^+ n^+},
\nn \\[\jot]
\piappr(b) \,\kett{j\mu n}
&= \bappr^+_{j\mu n} \kett{j^+ \mu^+ n^-}
 + \bappr^-_{j\mu n} \kett{j^- \mu^+ n^-},
\nn \\[\jot]
\piappr(a^*) \,\kett{j\mu n}
&= \tilde\aappr^+_{j\mu n} \kett{j^+ \mu^- n^-}
 + \tilde\aappr^-_{j\mu n} \kett{j^- \mu^- n^-},
\label{eq:appr-spin-repn}
\\[\jot]
\piappr(b^*) \,\kett{j\mu n}
&= \tilde\bappr^+_{j\mu n} \kett{j^+ \mu^- n^+}
 + \tilde\bappr^-_{j\mu n} \kett{j^- \mu^- n^+},
\nn
\end{align}
where
\begin{align}
\aappr^+_{j\mu n}
&:= \sqrt{1 - q^{2j+2\mu+2}}
\begin{pmatrix} \sqrt{1 - q^{2j+2n+3}} & 0 \\
0 & \sqrt{1 - q^{2j+2n+1}} \end{pmatrix},
\nn \\[\jot]
\aappr^-_{j\mu n}
&:= q^{2j+\mu+n+\half} \sqrt{1 - q^{2j-2\mu}}
\begin{pmatrix} q \sqrt{1 - q^{2j-2n+1}} & 0 \\
0 & \sqrt{1 - q^{2j-2n-1}} \end{pmatrix},
\nn \\[\jot]
\bappr^+_{j\mu n}
&:= q^{j+n-\half} \sqrt{1 - q^{2j+2\mu+2}}
\begin{pmatrix} q \sqrt{1 - q^{2j-2n+3}} & 0 \\
0 & \sqrt{1 - q^{2j-2n+1}} \end{pmatrix},
\label{eq:piappr-coeff}
\\[\jot]
\bappr^-_{j\mu n}
&:= - q^{j+\mu} \sqrt{1 - q^{2j-2\mu}}
\begin{pmatrix} \sqrt{1 - q^{2j+2n+1}} & 0 \\
0 & \sqrt{1 - q^{2j+2n-1}} \end{pmatrix},
\nn
\end{align}
and
\begin{equation}
\tilde{\aappr}^\pm_{j\mu n} = \aappr^\mp_{j^\pm \mu^- n^-},  \qquad
\tilde{\bappr}^\pm_{j\mu n} = \aappr^\mp_{j^\pm \mu^- n^+}.
\end{equation}
\end{prop}

\begin{proof}
First of all, we claim that the defining relations
\eqref{eq:suq2-relns} are preserved by $\piappr$ modulo the ideal
$\K_q$ of $\B(\H)$, that is,
$\piappr(b) \piappr(a) - q\,\piappr(a) \piappr(b) \in \K_q$,
and so on. Indeed, it can be verified by a direct but tedious check on
the spinor basis that
$\piappr(b) \piappr(a) - q\,\piappr(a) \piappr(b) = L_q^4 A$
where $A$ is a bounded operator; the same is true for each of the
other relations listed in~\eqref{eq:suq2-relns}.

It is well known, and easily checked from \eqref{eq:suq2-relns}, that
$\A$ is generated as a vector space by the products $a^k b^l b^{*m}$
and $b^l b^{*m} a^{*n}$, for $k,l,m,n \in \N$. We may thus define
$\piappr(x)$ for any $x \in \A$ by extending
\eqref{eq:appr-spin-repn} multiplicatively on such products, and then
extending further by linearity. With this convention, we conclude that
\begin{equation}
\piappr(xy) - \piappr(x) \piappr(y) \in \K_q
\sepword{for all}  x,y \in \A.
\label{eq:almost-repn}
\end{equation}
The defining formulas also entail that $\piappr(x)^* = \piappr(x^*)$
for each $x \in \A$.

If $\pi'(x) - \piappr(x) \in \K_q$ and
$\pi'(y) - \piappr(y) \in \K_q$, then
$$
\pi'(xy) - \piappr(x) \piappr(y)
= \pi'(x) \bigl( \pi'(y) - \piappr(y) \bigr)
+ \bigl( \pi'(x) - \piappr(x) \bigr) \piappr(y) \in \K_q,
$$
and therefore $\pi'(xy) - \piappr(xy)$ lies in $\K_q$ also; thus, it
suffices to verify this property in the cases $x = a,b$.

On comparing the coefficients \eqref{eq:piappr-coeff} with the 
corresponding ones of $\pi'(a)$ and $\pi'(b)$
from~\eqref{eq:spin-coeff}, we get, for instance,
\begin{subequations}
\label{eq:coeff-diff}
\begin{align}
\a^+_{j\mu n\up\up} - \aappr^+_{j\mu n\up\up}
&= \frac{q^{4j+4} \sqrt{1 - q^{2j+2\mu+2}} \sqrt{1 - q^{2j+2n+3}}}
{1 - q^{4j+4}}
= q^{4j+4} \,\a^+_{j\mu n\up\up},
\nn \\
\a^+_{j\mu n\dn\dn} - \aappr^+_{j\mu n\dn\dn}
&= \frac{q^{4j+2} \sqrt{1 - q^{2j+2\mu+2}} \sqrt{1 - q^{2j+2n+1}}}
{1 - q^{4j+2}}
= q^{4j+2} \,\a^+_{j\mu n\dn\dn}.
\label{eq:coeff-diff-plus}
\end{align}
and similarly,
\begin{equation}
\a^-_{j\mu n\up\up} - \aappr^-_{j\mu n\up\up}
= q^{4j+2} \,\a^-_{j\mu n\up\up},  \qquad
\a^-_{j\mu n\dn\dn} - \aappr^-_{j\mu n\dn\dn}
= q^{4j} \,\a^-_{j\mu n\dn\dn}.
\label{eq:coeff-diff-minus}
\end{equation}
\end{subequations}
We estimate the off-diagonal terms, using the inequalities
$q^{\pm\mu} \leq q^{-j}$, $q^{\pm n} \leq q^{-j-\half}$ and
$[N]^{-1} < q^{N-1}$:
\begin{align*}
|\a^+_{j\mu n\dn\up}|
&= q^{(\mu+n+\half)/2} \,
\frac{[j + \mu + 1]^\half \,[j - n + \half]^\half} {[2j+1]\,[2j+2]} 
\leq \frac{q^{-2j-2}}{[2j+1]\,[2j+2]} < q^{2j-1},
\\
|\a^-_{j\mu n\up\dn}|
&= q^{(\mu+n+\half)/2} \,
\frac{[j - \mu]^\half \,[j + n + \half]^\half} {[2j]\,[2j+1]} 
\leq \frac{q^{-2j-1}}{[2j]\,[2j+1]} < q^{2j-2}.
\end{align*}
On account of \eqref{eq:coeff-diff} and analogous relations for the
coefficients of $\piappr(b)$, we find that
\begin{align*}
\pi'(a) - \piappr(a) &\equiv T \pi'(a) T \quad \bmod \K_q,
\\ 
\pi'(b) - \piappr(b) &\equiv T \pi'(b) T \quad \bmod \K_q,
\end{align*}
where $T$ is the operator defined by
\begin{equation}	 
T\,\kett{j\mu n}
:= \begin{pmatrix} q^{2j+\sesq} & 0 \\ 0 & q^{2j+\half} \end{pmatrix}
\,\kett{j\mu n}
= \begin{pmatrix} q^\sesq & 0 \\ 0 & q^\half \end{pmatrix}
L_q^2 \,\kett{j\mu n}.
\label{eq:T-defn}
\end{equation}
Clearly, $T \in \K_q$, so that by boundedness of $\pi'(x)$ it follows
that $\pi'(x) - \piappr(x) \in \K_q$ for $x = a,b$.
\end{proof}

Using the conjugation operator $J$, we can also define an approximate
antirepresentation of $\A$ by
$\piiapprop(x) := J \piappr(x) J^{-1}$. It is immediate that
$\piiop(x) - \piiapprop(x) \in \K_q$, with $\piiop$ as
defined in Proposition~\ref{pr:piiop-eqvar}. Explicitly, we can write
\begin{align*}
\piiapprop(a) \,\kett{j\mu n}
&= \aappr^{\circ+}_{j\mu n} \kett{j^+ \mu^+ n^+}
 + \aappr^{\circ-}_{j\mu n} \kett{j^- \mu^+ n^+},
\\[\jot]
\piiapprop(b) \,\kett{j\mu n}
&= \bappr^{\circ+}_{j\mu n} \kett{j^+ \mu^+ n^-}
 + \bappr^{\circ-}_{j\mu n} \kett{j^- \mu^+ n^-},
\\[\jot] 
\piiapprop(a^*) \,\kett{j\mu n} 
&= \tilde\aappr^{\circ+}_{j\mu n} \kett{j^+ \mu^- n^-}
 + \tilde\aappr^{\circ-}_{j\mu n} \kett{j^- \mu^- n^-},
\\[\jot]
\piiapprop(b^*) \,\kett{j\mu n}
&= \tilde\bappr^{\circ+}_{j\mu n} \kett{j^+ \mu^- n^+}
 + \tilde\bappr^{\circ-}_{j\mu n} \kett{j^- \mu^- n^+},
\end{align*}
where 
$$
\aappr^{\circ\pm}_{j\mu n} = \tilde\aappr^\pm_{j,-\mu,- n}, \quad
\tilde\aappr^{\circ\pm}_{j\mu n} = \aappr^\pm_{j,-\mu,- n}, \quad
\bappr^{\circ\pm}_{j\mu n} = - \tilde\bappr^\pm_{j,-\mu,-n}, \quad
\tilde\bappr^{\circ\pm}_{j\mu n} = - \bappr^\pm_{j,-\mu,-n}.
$$

\vspace{6pt}

It turns out that the approximate representations $\piappr$ and
$\piiapprop$ almost commute, in the following sense.

\begin{prop} 
\label{pr:pi-comm}
For each $x,y \in \A$, the commutant
$[\piiapprop(x), \piappr(y)]$ lies in~$\K_q$.
\end{prop}

\begin{proof}
In view of our earlier remarks on the almost-multiplicativity of
$\piappr$, and thus also of~$\piiapprop$, it is enough to check this
for the cases $x, y = a, a^*, b, b^*$. We omit the detailed
calculation, which we have performed with a symbolic computer program.
In each case, the commutator $[\piiapprop(x), \piappr(y)]$ decomposes
as a direct sum of operators in the subspaces $W_j^\up$ and $W_j^\dn$
separately, in view of \eqref{eq:piappr-coeff}
and~\eqref{eq:J-formula}, and the explicit calculation shows that for
each pair of generators $x,y$, we obtain
$[\piiapprop(x), \piappr(y)] = L_q^2 A$ where $A$ is a bounded
operator.
\end{proof}

If we further impose the first-order condition up to compact operators
in the ideal $\K_q$, it turns out that this (almost) determines the
Dirac operator.

\begin{prop}
\label{pr:first-order}
Up to rescaling, adding constants, and adding elements of~$\K_q$,
there is only one operator $D$ of the form \eqref{eq:Dirac-eigen}
which satisfies the first order condition modulo~$\K_q$, that is,
each $[D, \piappr(y)]$ is bounded, and
\begin{equation}
[\piiapprop(x), [D, \piappr(y)]] \in \K_q
\sepword{for all} x,y \in \A.
\label{eq:approx-firstord}
\end{equation}
This operator $D$ has eigenvalues that are linear in~$j$. 
\end{prop}

\begin{proof}
Suppose first that $D$ is an equivariant selfadjoint operator of the
type considered in Section~\ref{sec:Dirac-oper}, with eigenvalues
linear in~$j$; that is, $D$ is determined by \eqref{eq:Dirac-eigen}
and~\eqref{eq:linear-evs}. Since each operator appearing in
\eqref{eq:approx-firstord} decomposes into a pair of operators on the 
``up'' and ``down'' spinor subspaces, it is clear that the nested 
commutators are independent of the parameters $c_2^\up$ and $c_2^\dn$;
and that $c_1^\up$ and $c_1^\dn$ are merely scale factors on both
subspaces. Again we take $x$ and $y$ to be generators: explicit
calculations show that in each case, 
$[\piiapprop(x), [D, \piappr(y)]] = L_q^2 B$ with $B$ a bounded
operator.

To prove the converse, assume only that $D$ satisfies the 
equivariance condition  \eqref{eq:Dirac-eigen}, and that
$[D, \piappr(a)]$ and $[D, \piappr(b)]$ are bounded.

We may decompose $\piappr(a) = \piappr(a)^+ + \piappr(a)^-$ according
to whether the index~$j$ in \eqref{eq:appr-spin-repn} is raised or
lowered; and similarly for $\piappr(b)$, $\piiapprop(a)$, and
$\piiapprop(b)$. Proposition~\ref{pr:pi-comm} shows that,
modulo~$\K_q$:
\begin{align*}
\piappr(a)^+ \piiapprop(a)^+ &\equiv \piiapprop(a)^+ \piappr(a)^+,
\\
\piappr(a)^- \piiapprop(a)^- &\equiv \piiapprop(a)^- \piappr(a)^-,
\\
\piappr(a)^+ \piiapprop(a)^- + \piappr(a)^- \piiapprop(a)^+
&\equiv \piiapprop(a)^+ \piappr(a)^- + \piiapprop(a)^- \piappr(a)^+ .
\end{align*}

By \eqref{eq:piappr-coeff}, the operators $\piappr(a)$ and
$\piappr(b)$, as well as $D$, are diagonal for the decomposition
$\H = \H^\up \oplus \H^\dn$. On the subspace $\H^\up$, we obtain
\begin{align}
& [[D, \piappr(a)], \piiapprop(a)] \,\ket{j\mu n\up}
\nn \\
&= \bigl( D \piappr(a) \piiapprop(a)
+ \piiapprop(a) \piappr(a) D - \piappr(a) D \piiapprop(a)
- \piiapprop(a) D \piappr(a) \bigr) \,\ket{j\mu n\up}
\nn \\
&= \Bigl(
(d_{j+1}^\up + d_j^\up - 2d_{j^+}^\up) \,\piappr(a)^+ \piiapprop(a)^+
+ (d_{j-1}^\up + d_j^\up - 2d_{j^-}^\up) \,\piappr(a)^- \piiapprop(a)^-
\nn \\
&\qquad +
2 d_j^\up (\piappr(a)^+ \piiapprop(a)^- + \piappr(a)^- \piiapprop(a)^+)
- d_{j^+}^\up (\piappr(a)^- \piiapprop(a)^+
+ \piiapprop(a)^- \piappr(a)^+ )
\nn \\
&\qquad - d_{j^-}^\up
(\piappr(a)^+ \piiapprop(a)^- + \piiapprop(a)^+ \piappr(a)^-) 
+ R \Bigr) \,\ket{j\mu n\up},
\label{eq:D-bicomm}
\end{align}
where $R \in \K_q$. On the subspace $\H^\dn$, we get the precisely
analogous expression with the arrows reversed.

In order that the expression on the right hand side
of~\eqref{eq:D-bicomm} come from an element of $\K_q$ applied to
$\ket{j\mu n\up}$, and likewise for $\ket{j\mu n\dn}$, it is
necessary and sufficient that the scalars
\begin{equation}
w_j^\up := d_{j+1}^\up + d_j^\up - 2d_{j^+}^\up,  \qquad
w_j^\dn := d_{j+1}^\dn + d_j^\dn - 2d_{j^+}^\dn
\label{eq:almost-linear}
\end{equation}
satisfy $w_j^\up = O(q^j)$ and $w_j^\dn = O(q^j)$ as $j \to \infty$.

In the particular case where $w_j^\up = 0$ and $w_j^\dn = 0$ for
all~$j$, \eqref{eq:almost-linear} gives elementary recurrence
relations for $d_j^\up$ and $d_j^\dn$, whose solutions are precisely
the expressions \eqref{eq:linear-evs} that are linear in~$j$, namely,
$$
d_j^\up = c_1^\up j + c_2^\up, \qquad
d_j^\dn = c_1^\dn j + c_2^\dn.
$$
The general case gives a pair of perturbed recurrence relations, that
may be treated by generating-function methods \cite{GrahamKP}; their
solutions differ from the linear case by terms that are $O(q^j)$ as
$j \to \infty$. Thus, the corresponding operator $D$ differs from one
whose eigenvalues are linear in~$j$ by an element of~$\K_q$.
\end{proof}

We finish by summarizing the implications of the above Propositions
\ref{pr:appr-rep}, \ref{pr:pi-comm} and~\ref{pr:first-order} for the
spectral triple $(\A(SU_q(2)),\H,D,J)$, where $\A(SU_q(2))$
acts on $\H$ via the spinor representation~$\pi'$. 

The representations $\pi'$ and $\piiop$ do not commute, since
the conjugation operator~$J$ differs from the Tomita conjugation
for~$\pi'$. However, we do obtain commutation ``up to infinitesimals''; since
$[\piiop(x), \pi'(y)] \equiv [\piiapprop(x), \piappr(y)] \bmod \K_q$,
Proposition~\ref{pr:pi-comm} entails the analogous result for the
exact representations:
$$
[\piiop(x), \pi'(y)] \in \K_q  \sepword{for all}  x,y \in \A.
$$

To examine the first-order property, we note first if $x,y \in \A$ and
$[D, \pi'(y) - \piappr(y)]$ lies in~$\K_q$, then
\begin{align} 
[\piiop(x), [D,\pi'(y)]]
&= \bigl[ \piiapprop(x) + (\piiop(x) - \piiapprop(x)), \
[D,\piappr(y) + (\pi'(y) - \piappr(y))] \bigr] 
\nn \\
&\equiv [\piiapprop(x), [D, \piappr(y)]] \equiv 0 \bmod \K_q.
\label{eq:almost-firstord}
\end{align}
Since $D$ commutes with the positive operator $T$ defined in
\eqref{eq:T-defn}, we find in the case of a generator $y = a$, $a^*$,
$b$ or~$b^*$, that
$$
[D, \pi'(y) - \piappr(y)] = [D, T \pi'(y) T] = T [D,\pi'(y)] T,
$$
which lies in $\K_q$ since $[D, \pi'(y)]$ is bounded, by
Proposition~\ref{pr:D-comm}. Thus, $[D, \piappr(y)]$ is bounded, too
--as required by Proposition~\ref{pr:first-order}. The general case
of $[D, \pi'(y) - \piappr(y)] \in \K_q$ then follows
from~\eqref{eq:almost-repn}. Thus \eqref{eq:almost-firstord} holds for
general $x,y \in \A$. Combining that with
Proposition~\ref{pr:first-order}, we arrive at the following
characterization of our spectral triple over $\A(SU_q(2))$.

\begin{thm}
\label{th:grs}
The real spectral triple $(\A(SU_q(2)), \H, D, J)$ defined here,
with $\A(SU_q(2))$ acting on $\H$ via the spinor representation
$\pi'$, satisfies both the commutant property and the first order
condition up to infinitesimals:
$$
\begin{aligned}
{}[\piiop(x), \pi'(y)] \in \K_q,
\\[\jot]
[\piiop(x), [D,\pi'(y)]] \in \K_q,
\end{aligned}
\quad \sepword{for all }  x,y \in \A(SU_q(2)).
$$
\end{thm}

In \cite{Hawkins} it was argued that there are obstructions to the
construction of ``deformed spectral triples'' satisfying a type of
first order condition for the Dirac operator. Theorem~\ref{th:grs}
above shows a way to overcome these obstructions.

\subsection*{Acknowledgments}

We thank 
Alain Connes and Mariusz Wodzicki for helpful discussions.
AS and JCV thank Alain Connes for the opportunity to visit the
Institut des Hautes \'Etudes Scientifiques. AS also thanks the
Institut Henri Poincar\'e and SISSA for hospitality. LD acknowledges a
partial support by the EU Project INTAS 00-257. Support from the
Vicerrector\'{\i}a de Investigaci\'on of the Universidad de Costa Rica
is also acknowledged. JCV is grateful to the DAAD and to Florian
Scheck for a warm welcome to the Johannes-Gutenberg Universit\"at,
Mainz, during the course of this work.

\end{document}